\documentclass[12pt]{article}

\usepackage{authblk}

\usepackage{amsmath,amssymb,amsthm}
\usepackage{booktabs}
\usepackage{multirow}
\usepackage{makecell}
\usepackage{enumerate}
\usepackage{graphicx}
\graphicspath{{figs/}}
\usepackage{url}
\usepackage{color,subcaption}
\usepackage{geometry}
\usepackage{algorithm}
\usepackage{algpseudocode}
\usepackage{caption}
\usepackage[counterclockwise]{rotating}
\usepackage[T1]{fontenc}
\usepackage{setspace}
\usepackage{thmtools}
\usepackage{thm-restate}
\usepackage[bookmarksopen=false]{hyperref}
\usepackage[flushleft]{threeparttable}
\usepackage[round]{natbib}
\usepackage{bm} 

\usepackage{subcaption} 
\usepackage{cleveref}
\usepackage{diagbox} 

\usepackage{tikz}
\hypersetup{pdfborder = {0 0 0},colorlinks=true,linkcolor=blue,citecolor=blue}
\newtheorem{assumption}{Assumption}
\newtheorem{theorem}{Theorem}

\newtheorem{lemma}{Lemma}
\newtheorem{corollary}{Corollary}
\theoremstyle{definition}

\newcommand{\amatch}{\hat{A}_{match}}
\newcommand{\amatchbc}{\hat{A}^{bc}_{match}}
\newcommand{\nmeani}{\frac{1}{n}\sum_{i=1}^n}

\newcommand{\Yhat}{\hat{Y}}
\newcommand{\Ytilde}{\tilde{Y}}
\newcommand{\KMpi}{K_M(\pi,i)}
\newcommand{\KMi}{K_M(i)}
\newcommand{\JMi}{\mathcal{J}_M(i)}

\newcommand{\muXOne}{\bm{\mu}\mathbf{(X \text{,} 1)}}
\newcommand{\muXZero}{\bm{\mu}\mathbf{(X \text{,} 0)}}

\allowdisplaybreaks[4]

\title{Matching-Based Policy Learning}

\author{Xuqiao Li}
\affil{School of Mathematics, Sun Yat-sen University, Guangzhou, China, 510275\\
	Email: \href{mailto: lixq87@mail2.sysu.edu.cn}{\textcolor{black}{lixq87@mail2.sysu.edu.cn}}}

\author{Ying Yan\thanks{Corresponding Author.}}
\affil{School of Mathematics, Sun Yat-sen University, Guangzhou, China, 510275\\
	Email: \href{mailto: yanying7@mail.sysu.edu.cn}{\textcolor{black}{yanying7@mail.sysu.edu.cn}}}
\date{}

\begin{document}
	\maketitle

\thispagestyle{empty}

\begin{abstract}
\noindent The beneficial effects of  treatments vary across individuals in most  studies. Treatment heterogeneity  motivates practitioners to search for the optimal policy  based on personal characteristics. A long-standing common practice in policy learning has been estimating and  maximizing the  value function using weighting techniques. Matching is  widely used  in many applied disciplines to infer causal effects, which is   intuitively appealing  because the observed covariates   are directly  balanced across different treatment groups. Nevertheless, matching is rarely explored in policy learning. In this work,
 we propose a matching-based policy learning framework. We adapt  standard and bias-corrected  matching methods to estimate an alternative form of the value function: the  advantage function, which can be interpreted as  the expected
improvement achieved by implementing a given policy  compared to the equiprobable random policy. We then learn the optimal policy over a restricted policy class by maximizing the matching estimator of the advantage function.  We derive a non-asymptotic high probability bound for the regret of the learned optimal policy. Moreover, we show that  the learned  policy  is almost rate-optimal. The competitive finite sample performance of the proposed method compared to weighting-based and outcome modeling-based learning methods is demonstrated in extensive simulation studies and a real data application.
\end{abstract}


\section{Introduction}
The beneficial effects of  treatments or interventions  vary across individuals in most studies. For example, a training program that increases employment rate for some workers who would otherwise have remained unemployed, but others may remain  unemployed despite intervention. Recently, data-driven personalized decision-making has received increased attention in randomized controlled trials and observational studies in various fields. In healthcare, clinicians utilize  patients' characteristics and medical history to tailor diabetes management to improve health outcomes \citep{Bertsimas2017}. In public policy, leveraging individual features  lead to  improved allocation of social services and thus enhance  overall efficiency \citep{KubeDasFowler2019}. In product recommendations, the online platforms  effectively  increase customer retention by recommending the products that match customers' personal preferences \citep{bastani2022learning}. A common underlying feature of these applications is the presence of treatment heterogeneity, which is a blessing that brings opportunities to improve the expected outcome in the target population. In view of  substantially increasing volumes of personal data collected by advanced technologies nowadays, it is  imperative to  exploit  heterogeneity of treatment effects efficiently and search for the optimal policy that allocates treatments based on individual characteristics.

Statistical methods for learning the optimal policy have been extensively developed for more than a decade. They can be roughly divided into two categories: model-based methods and direct-search methods. Q-learning \citep{qian2011performance}   and A-learning \citep{shi2018HD_Alearning} are model-based methods, which first fit the mean outcome given covariates and treatment or posit a model for the contrast function, and then determine the optimal policy based on the corresponding estimate. The performance of these methods critically relies  on correct  specification of the outcome model or the contrast function. Besides Q-learning, a long-standing common practice in policy learning has been estimating and  maximizing the  value function using weighting techniques. Outcome weighted learning \citep[O-learning,][]{zhao2012estimating} is a pioneering work of the direct-search methods that nonparametrically learns the optimal policy by utilizing  inverse probability weighting  (IPW) to estimate the value function in randomized clinical trials. In observational studies, the unknown propensity score in O-learning needs to be estimated using logistic regression or machine learning approaches \citep{lee2010improving}. One recognized potential limitation of  weighting-based approaches is the numerical instability caused by extreme values of the estimated propensity scores or model misspecification.

To enhance efficiency and robustness of IPW-based  learning methods, augmented inverse probability weighting  (AIPW) has been introduced  to policy learning \citep[e.g.,][]{zhangb2012stat,zhangb2012biom,zhao2019efficient,athey2021policy,pan2021improved}. \citet{zhangb2012biom} proposed to use AIPW  to estimate the value function and showed that the estimation procedure is doubly robust in the sense that the estimator is consistent if either the  parametric propensity score model or the parametric outcome mean model is correctly posited. However, the performance may be unsatisfactory  if both models are misspecified \citep{kangschafer2007}. \citet{athey2021policy} introduced AIPW-based policy learning using semiparametric theory, which leverages the cross-fitting technique \citep{Chernozhukov2018DoubleML} and flexible machine learning methods to  nonparametrically estimate the nuisance functions: the propensity score  and  the  conditional mean outcome. They proposed the policy tree algorithm to search for the optimal policy over a policy class. The learned policy is shown to achieve minimax optimal regret if the mean-squared errors of the nuisance functions satisfy certain order conditions of convergence rate.
Despite the promising theoretical guarantee, the performance of AIPW-based learning  deteriorates when extreme weights are present  or  when the sample size is  small or moderate.  Similar discussions can be found  in the literature \citep[e.g.,][]{zhao2019efficient,wu2022matching}. When the underlying conditional average treatment effect is not continuous at some points in the  covariate space, \citet{athey2021policy} observed the interesting  ``phase transition'' phenomenon, in which  AIPW-based policy learning exhibited enormous  variation in the estimated value function with the learned optimal policy when the sample size is moderate, whereas its performance is satisfactory for either small or large sample size.

Matching is  widely used  in many applied disciplines to infer causal effects \citep{imbens_rubin_2015, rosenbaum2020modern}. The goal of matching is to produce covariate balance  such that the distributions of observed covariates in different treatment groups are approximately the same as they would be in a successful randomized experiment, thereby  mimicking randomized experiments to the best extent with observational data \citep{Stuart2010, lin2023estimation, YangZhang2023}. Therefore,  matching has clearer  intuitive appeal to practitioners compared to weighting techniques such as IPW and AIPW. Moreover, it is more robust to extreme values of the propensity score function  in many circumstances \citep[e.g.,][]{wu2022matching,YangZhang2023}. Nevertheless, matching is rarely explored in policy learning except the work by \cite{wu2020matched}.

In this paper, we develop  a novel matching-based framework for policy learning in observational studies with favorable finite sample performance. The key idea behind the proposed matching-based policy learning (MB-learning) originates from an alternative form of the value function based on the pair of potential outcomes (referred to as advantage function in Section \ref{section methodology}), which can be interpreted as   the  improvement on the expected outcome after implementing a given policy compared to a random policy with equal probabilities. We nonparametrically impute the missing potential outcomes using the response data within the matched set to develop a  matching estimator of  the advantage function, where we utilize  nearest-neighbor (NN) matching on covariates with replacement \citep{abadie2006large}. One recognized limitation of  matching estimator is its susceptibility to the curse of dimensionality, and it is not  $\sqrt{n}$-consistent due to a conditional  bias term incurred by sampling with replacement during matching.  To mitigate this issue, we employ the bias correction technique in \citet{abadie2011bias} and develop an improved matching estimator of the advantage function with $\sqrt{n}$-consistency.  Leveraging the formulation of the matching estimator, we  recast the optimization problem of  policy learning into a weighted classification problem and  derive the optimal policy using the policy tree algorithm \citep{athey2021policy}.

This paper makes the following major contributions to the existing literature.
First, different from the (A)IPW-based approaches, MB-learning is a  robust method for learning the optimal policy in observational studies. By utilizing a bias-corrected matching-based  estimator, MB-learning mitigates the potential variability arising from propensity score model misspecification or extreme weights.
Second, we establish theoretical guarantee for the optimal policy learned by the matching-based estimator, where a non-asymptotic high probability bound for the value function difference between the best-in-class policy and the learned policy is derived. Moreover, we demonstrate the superiority of MB-learning both theoretically and empirically.  \citet{athey2021policy} established that   AIPW-based policy learning with cross-fitting achieves rate-optimality in terms of minimax regret with the order of  $\sqrt{\text{VC}(\Pi)/n}$, where $\text{VC}(\Pi)$ describes the complexity of a given policy class $\Pi$ which will be defined formally in the following sections. In comparison, we prove that the regret bound of MB-learning is almost of order $\sqrt{\text{VC}(\Pi)/n}$ for any fixed  size of the matched set. Although MB-learning is not exactly rate-optimal,  it exhibits more robust finite sample performance in simulation studies and real data analysis especially when the sample size is not large or when the policy allocations are extremely imbalanced. Moreover, the undesirable ``phase transition'' phenomenon encountered when implementing  AIPW-based policy learning with cross-fitting does not occur when the proposed MB-learning is in use.

To the best of  our knowledge, matching was rarely utilized for  policy learning except  matched learning \citep[M-learning]{wu2020matched}, where  a general matched pair-based objective function is formulated and  shown to be  consistent of the value function for some special case and  the optimal policy is learned by weighted support vector machine. We highlight several distinctive features of MB-learning compared to M-learning. First, rather than  estimating the value function or its equivalent counterpart directly, the first step of M-learning is to construct a  complicated matched pair-based objective function and then prove its consistency in theory. In contrast, MB-learning is more intuitive  because it \textit{directly} targets the advantage function and uses various  matching techniques to obtain reasonable estimation. Second, we acknowledge and tackle the conditional bias due to matching, correlation due to sampling with replacement, and we  fix the number of matches in the   establishment of the asymptotic property and the non-asymptotic regret bound for MB-learning. In  M-learning, there is a lack of discussions on correlation due to sampling with replacement, and  the number of matches is forced to go to infinity with the sample size in  theory. Third, the policy derived from MB-learning possesses appealing interpretability  for practitioners because policy tree \citep{athey2021policy} is used in the optimization procedure, whereas  M-learning with  Gaussian kernels may be less interpretable. Finally, the  bias introduced by NN matching may substantially  deteriorate the performance of M-learning when moderate or high-dimensional covariates are present, whereas potential improvement using bias-corrected  techniques remains unexplored in M-learning. In comparison, we mitigate this issue using the bias-corrected matching estimator \citep{abadie2011bias} with Lasso in MB-learning.

The rest of this paper is organized as follows. In Section \ref{section methodology}, we develop the methodology of MB-learning. We first formulate the policy learning framework,   introduce the advantage function, and describe AIPW-based  policy learning  developed by \citet{athey2021policy}. Then, we propose a  matching estimator of the advantage function. The estimator is further improved by a bias correction technique. In Section \ref{section implementation}, we formulate MB-learning  from a weighted classification perspective, and derive the optimal policy via policy tree using the bias-corrected matching estimator of the advantage function. We establish rigorous asymptotic and non-asymptotic theoretical guarantees for MB-learning. The empirical performance of MB-learning is demonstrated through extensive simulation studies in Section \ref{section simulation}. An application to National Supported Work (NSW) Program is presented in Section \ref{section application}. We conclude this article with discussions in Section \ref{section discussion}.  The proofs of the theories and additional numerical results are relegated to the Supplementary Material.

\section{Methodology}\label{section methodology}
\subsection{Notations and Preliminaries}\label{section notation}
Let $X \in \mathcal{X}$ be a vector of $p$-dimensional pre-treatment covariates, where $\mathcal{X} \subset \mathbb{R}^p$ is a covariate space, $W \in \{0,1\}$  is a binary policy or treatment indicator, and $Y$ is an outcome of interest. Without loss of generality, we assume that a larger $Y$ is more desirable.  Under the potential outcome framework \citep{imbens_rubin_2015}, we define $Y(w)$ to be the potential outcome if the unit had received the treatment $W=w$. We adopt the standard assumptions in causal inference: (i) consistency: $Y=Y(W)$; (ii) unconfoundedness: $\{Y(1),Y(0)\} \perp  W \mid X$; and (iii) overlap: for any $x \in \mathcal{X}$, the propensity score $e(x)=\Pr(W=1 \mid X=x)$ satisfies  $c < e(x) < 1-c $ for some positive constant $c$. Additionally, let $\mu(x,w)=E\left[Y|X=x, W=w\right]$, $\sigma^2(x,w)=\text{Var}\left[Y|X=x, W=w\right]$, and $\epsilon_i = Y_i-\mu(X_i,W_i)$. Under the assumptions of consistency and unconfoundedness, we have $\mu(x,w)=E\left[Y(w)|X=x\right]$ and $\sigma^2(x,w)=\text{Var}\left[Y(w) \mid X=x\right]$.  Let $\tau=E[Y(1)-Y(0)]$ and $\tau(X)=E\left[Y(1)-Y(0)\mid X\right]$ represent the average treatment effect (ATE) and conditional average treatment effect (CATE), respectively. Furthermore, we assume that the observed data $\left\{(X_i,W_i,Y_i)\right\}_{i=1}^n$ are $n$ independent and identically distributed (i.i.d.) copies of $(X,W,Y)$.

We consider the problem of learning the optimal policy in observational studies, where  we define that any policy $\pi$ is   a deterministic  decision function mapping from  the covariate space $\mathcal{X}$ to the treatment set $\{0,1\}$. In practice, it is  unrealistic to assign a complicated  policy to the target population without restriction. One should attach problem-specific constraints to the  policy, such as political restriction \citep{kitagawa2018Econometrica}, fairness \citep{athey2021policy} or interpretability \citep{zhangb2012stat}. Therefore, we seek to learn the optimal policy within a restricted policy class $\Pi$.  

The value function $V(\pi)$ for the policy $\pi$ \citep[e.g.,][]{zhao2012estimating} is the expected outcome after implementing the treatment allocation by $\pi$, defined as
$V(\pi)=E\left[Y(\pi(X))\right]$.
The optimal policy $\pi^*$ within the class $\Pi$ maximizes the value function $V(\pi)$, i.e.,
$\pi^*=\arg\max_{\pi \in \Pi} V(\pi)$. The advantage function
\begin{equation*}\label{advantage function}
A(\pi)=E\left[\left(2\pi(X)-1\right)\tau(X)\right]
\end{equation*}
serves as an alternative form of the value function. Because $V(\pi)=E[Y(0)+(Y(1)-Y(0))\pi(X)]$, it follows that 
$$A(\pi)=2\left\{V(\pi)- \frac{E[Y(1)]+E[Y(0)]}{2}\right\}.$$ 
Hence, the advantage function $A(\pi)$ can be interpreted as twice  the improvement achieved by implementing policy $\pi$ relative to a randomly assigned policy  with equal probability \citep{mo2021JASADistributionalRobustITR}. That is,  it is a measure of evaluating the gain of implementing a given policy compared to the equiprobable random policy. The policy that maximizes $A(\pi)$ is equivalent to that maximizes $V(\pi)$. Throughout this paper, we focus on policy learning based on $A(\pi)$.

By leveraging  the idea of semiparametric efficient estimation \citep{robins1994estimation},  AIPW-based estimators of $V(\pi)$ or  $A(\pi)$ were developed in the  literature, and   the optimal policy was determined based on these estimators \citep[e.g.,][]{zhangb2012stat,athey2021policy}. The AIPW-based estimator of the advantage function proposed by \citet{athey2021policy} is 
$$\hat{A}_{aipw}(\pi) = \frac{1}{n} \sum_{i=1}^n \left\{2\pi(X_i)-1\right\} \Gamma_i,$$ where
$$
 \Gamma_i=\mu(X_i,1)-\mu(X_i,0)+ \frac{W_i-e(X_i)}{e(X_i)\left\{ 1-e(X_i)\right\} } \left\{Y_i - \mu(X_i,W_i)\right\}.
$$
This AIPW-based estimator of $A(\pi)$ is  efficient and doubly robust  against model misspecification in observational studies.
\citet{athey2021policy} employed black-box machine learning methods and the cross-fitting technique \citep{Chernozhukov2018DoubleML} to  nonparametrically estimate the nuisance functions: the propensity score $e(X_i)$  and  the  conditional mean outcome $\mu(X_i,W_i)$. The optimal policy is learned by maximizing the $\hat{A}_{aipw}(\pi)$ with the plugged-in estimated nuisance functions using the policy tree algorithm. 
 
Nevertheless, the procedure of estimating the nuisance functions may introduce instability particularly when the estimated propensity score is extreme. In addition,   machine learning methods may perform unsatisfactorily with a relatively small or moderate sample size. Weighting does not directly ensure covariate balance across treatment groups and is less intuitively appealing to practitioners as matching.
Therefore, it is imperative to develop a  sensible approach for learning the optimal policy  with ensured covariate balance property and desirable efficiency and robustness performance in finite samples. Learning policies using matching techniques is a natural choice.

\subsection{Matching-Based Advantage Function}\label{section Matching-based Advantage Function}

We rewrite  the advantage function in the following expression:
\begin{eqnarray*}
	A(\pi) &=& E\left[\left(2\pi(X)-1\right)\left\{Y(1)-Y(0)\right\}\right]\\
	      &=& E\left[\left\{2\pi(X)-1\right\}\operatorname{sign}\left\{Y(1)-Y(0)\right\}\left|Y(1)-Y(0)\right|\right].
\end{eqnarray*}
The second equality indicates that, in order to maximize  $A(\pi)$ one should assign each unit to the treatment group in which its  potential outcome is larger. Meanwhile, a more significant discrepancy between the pair of potential outcomes contributes to a greater influence in the treatment allocation. Building on this observation, we seek to develop an estimator of $A(\pi)$ based on the pair of potential outcomes. Half of the potential outcomes are missing, and we utilize matching techniques to impute the missing  outcomes.

We leverage the matching estimator in \citet{abadie2006large}, where  NN matching with replacement is in use. Specifically, for any $x \in \mathcal{X}$ and some positive definite matrix $V$, let $\Vert x \Vert = (x^\top V x)^{1/2}$ denote a vector norm.  We  define the matched set $\mathcal{J}_M(i)$ as the set that contains the indices of the first $M$ matches for unit $i$:
\begin{equation*}
	\mathcal{J}_M(i)=\left\{j:W_j=1-W_i, \sum_{k: W_k=1-W_i}I\{ \Vert X_k-X_i \Vert\leq \Vert X_j-X_i\Vert\} \leq M\right\}.
\end{equation*}
Throughout, we assume there are no ties in matching, which is satisfied if $X$ is continuously distributed \citep{abadie2006large}. Tackling ties will be theoretically more demanding and it is an interesting direction of future work. In this paper, we consider the  commonly used Mahalanobis metric:  $\Vert x \Vert = (x^\top V_{\text{maha}} x)^{1/2}$, where
$V_{\text{maha}}= \left\{ n^{-1} \sum_{i=1}^n (X_i - \bar{X})(X_i - \bar{X})^\top \right\}^{-1}$ and $\bar{X} = n^{-1} \sum_{i=1}^n X_i$. 
See \citet{abadie2006large}  for discussion on other metrics. The matching estimator of potential outcomes is given by
\begin{equation*}
	\hat{Y}_i(0)=\left\{ \begin{aligned}
	&Y_i, &\text{ if } W_i = 0,\\
	&\frac{1}{M} \underset{j \in \mathcal{J}_M(i)}{\sum} Y_j, &\text{ if } W_i = 1,
\end{aligned} \right.
\quad \text{ and } \quad
\hat{Y}_i(1)=\left\{ \begin{aligned}
	&\frac{1}{M} \underset{j \in \mathcal{J}_M(i)}{\sum} Y_j, &\text{ if }
	W_i = 0 ,\\ &Y_i, &\text{ if } W_i = 1.
\end{aligned} \right.
\end{equation*}
For each unit $i$, the pair of matching estimators $\{\hat{Y}_i(1),\hat{Y}_i(0)\}$ mimics the pair of potential outcomes $\{Y_i(1),Y_i(0)\}$. Finally, we propose the matching-based advantage function:
\begin{equation}\label{Afun_match}
	\hat{A}_{match}(\pi)= \frac{1}{n}\sum^n_{i=1} \left\{2\pi(X_i)-1\right\} \left\{\hat{Y}_i(1)-\hat{Y}_i(0)\right\}.
\end{equation}

In the following, we  establish the consistency and asymptotic normality of the matching-based estimator $\amatch(\pi)$, followed by discussion of its limitation and the proposal of improvement. Let $M_i=(M_{i1},\ldots,M_{in})^\top$ denote the $n$-dimensional vector indicating the units in $\mathcal{J}_M(i)$, where $M_{ij}=I\left\{j \in \mathcal{J}_M(i)\right\}$. That is, $M_{ij}=1$ if unit $j \in \mathcal{J}_M(i)$, and $M_{ij}=0$ otherwise.  For convenience of further analysis,  we rewrite $\hat{A}_{match}(\pi)$ as a linear combination of $\{Y_i\}_{i=1}^n$ using elementary algebra.
\begin{lemma}\label{linear_combination}
The matching-based advantage function $\hat{A}_{match}(\pi)$ has the following expression:
\begin{equation*}\label{eq_linear_combination}
	\hat{A}_{match}(\pi)= \frac{1}{n}\sum^n_{i=1} (2W_i-1)\left[\{2\pi(X_i)-1\}+\frac{K_M(\pi,i)}{M} \right]Y_i,
\end{equation*}
where $K_M(\pi,i)= \sum^n_{j=1}\{2\pi(X_j)-1\}M_{ji}=\sum_{j:M_{ji}=1}{2\pi(X_j)-1}$.
\end{lemma}
It is noteworthy that $\{K_M(\pi,i)\}_{i=1}^n$ is the key component for theoretical  developments. Recall that the matching estimator for the ATE $\tau$ in \citet{abadie2006large} is
\[\hat{\tau} =\frac{1}{n}\sum^n_{i=1}\left\{\hat{Y}_i(1)-\hat{Y}_i(0)\right\}= \frac{1}{n}\sum^n_{i=1} (2W_i-1)\left\{1+\frac{K_M(i)}{M}\right\}Y_i, \] where $\KMi= \sum_{j=1}^n M_{ji}$ is the number of times unit $i$ has been matched. The ATE estimator $\hat{\tau}$ and the advantage function estimator $\hat{A}_{match}(\pi)$ have similar expressions. In particular, when $\pi(x) \equiv 1$, it follows that $\hat{A}_{match}(\pi)=\hat{\tau}$ and $\KMpi=\KMi$. Therefore, the matching estimator $\hat{\tau}$ is a special case of our proposed  estimator with the ``treat-all'' policy. In the following, we show that the properties of $\amatch(\pi)$ resemble those of $\hat{\tau}$.

\begin{assumption}\label{assumption_support}
The  $p$-dimensional vector $X$ is continuously distributed on a convex and compact support $\mathcal{X}_0 \subset \mathbb{R}^p$. Its density is upper bounded.  Moreover, its density is bounded away from zero on $\mathcal{X}_0$.
\end{assumption}

Assumption 1 is a standard regularity condition for the matching estimator of average causal effect \citep{abadie2006large}. Discrete covariates with finite supports can be addressed by analyzing the estimation of the advantage function within subsamples. For high-dimensional and  continuous $X$, the assumption that its density is bounded away from zero may be too stringent \citep{otsu2017bootstrap}. Future work on relaxing this assumption is warranted.

\begin{lemma}\label{lemma_KMpi_order}
  Under  Assumption 1, for any  $\pi \in \Pi$, $K_M(\pi,i)=O_p(1)$ and $E[K_M(\pi,i)^q]$ is bounded uniformly   for any $n$ and $q>0$.
\end{lemma}
This lemma says that  $\KMpi$ has the same order as $\KMi$ under   Assumption 1. Analyzing the asymptotic properties of $\amatch(\pi)$ is challenging due to the dependence among the terms $\{\KMpi\}^n_{i=1}$. To alleviate  the difficulty,  we decompose $\amatch(\pi)$ into several components for further analysis as in \citet{abadie2006large}.
\begin{lemma}\label{lemma_Amatch_decomposition}
The matching-based advantage function $\amatch(\pi)$ has the following decomposition:
\begin{equation*}
	\hat{A}_{match}(\pi)=\bar{A}(\pi)+E_M(\pi)+B_M(\pi),
\end{equation*}
where
\begin{eqnarray} \bar{A}(\pi)&=&\frac{1}{n}\sum^n_{i=1}\{2\pi(X_i)-1\}\{\mu(X_i,1)-\mu(X_i,0)\},\nonumber\\
		E_M(\pi)&=&\frac{1}{n}\sum^n_{i=1} (2W_i-1)\left[\{2\pi(X_i)-1\}+\frac{K_M(\pi,i)}{M} \right]\epsilon_i,\nonumber\\
	B_M(\pi)&=&\frac{1}{n}\sum^n_{i=1} (2W_i-1)\{2\pi(X_i)-1\}\left[\frac{1}{M}\underset{j \in \mathcal{J}_M(i)}{\sum}\{\mu(X_i,1-W_i)-\mu(X_j,1-W_i)\}\right].\label{eq_BMpi}
\end{eqnarray}
\end{lemma}

Let $\mathbf{X}$ be an $n \times p$ matrix where the $i$th row is $X^{\top}_i$. Let  $\mathbf{W}$ be the $n \times 1$ vector with the $i$th element equal to $W_i$. It follows that $B_M(\pi)=E[\amatch(\pi)-\bar{A}(\pi) \mid \mathbf{X}, \mathbf{W}]$ after some algebra, referred to as the \textit{conditional bias} term relative to $\bar{A}(\pi)$. Utilizing the decomposition in Lemma \ref{lemma_Amatch_decomposition}, we proceed to derive the  large sample properties of $\amatch(\pi)$ by analyzing its three components respectively. First,
we  derive the order of $B_M(\pi)$ in the following lemma.
\begin{lemma}\label{order_conditional_bias}
	Under    Assumption 1, suppose that $\mu(x,1)$ and $\mu(x,0)$ are Lipschitz in $\mathcal{X}_0$, then for any  $\pi \in \Pi$, $B_M(\pi)=O_p\left(n^{-1/p}\right)$.
\end{lemma}

Lemma \ref{order_conditional_bias} suggests that  the conditional bias $B_M(\pi)$ is asymptotically negligible, but the convergence rate is quite slow when the covariate dimension $p$ is moderate or large. The conditional variance of $E_M(\pi)$ given $\mathbf{X}$ and $\mathbf{W}$ is: 
$$\text{Var}\left[E_M(\pi) \mid \mathbf{X},\mathbf{W}\right]=\frac{1}{n^2}\sum^n_{i=1} \left[\{2\pi(X_i)-1\}+\frac{K_M(\pi,i)}{M}\right]^2 \sigma^2(X_i,W_i).$$ 
Define $V^E(\pi)=n \text{Var}\left[E_M(\pi) \mid \mathbf{X},\mathbf{W}\right]$.  

\begin{assumption}Assume the following conditions:

(i) For $w=0,1$, $\mu(x,w)$ and $\sigma^2(x,w)$ are Lipschitz in $\mathcal{X}_0$.

(ii) For $w=0,1$, E$\left[Y^4 \mid X=x,W=w \right]$ is upper bounded uniformly in $\mathcal{X}_0$.

(iii) For any $\pi \in \Pi$, $V^E(\pi)$ is bounded away from zero for any   $(\mathbf{X},\mathbf{W})\in\mathcal{X}_0^n  \times \{0,1\}^n$.
\end{assumption}

Assumption 2(i) is a mild condition  on  the smoothness and boundedness of the conditional mean and the conditional variance. Assumption 2(ii) is a moment condition for the central limit theorem. Assumption 2(iii)  is analogous to the Assumption 4(iii) in \citet{abadie2006large}, which guarantees the  non-singularity of  $V^E (\pi)$ uniformly in $\Pi$.

 Let $V^{\tau(X)}(\pi)=E[\{2\pi(X)-1\}\tau(X)-A(\pi)]^2$. We establish the asymptotic properties for $\hat{A}_{match}(\pi)$ in the following theorem.

\begin{theorem}\label{consistency and normality without bias correction}
(i)  Under   Assumptions 1 and 2(i), then for any  $\pi \in \Pi$, as $n\rightarrow\infty$, $$\hat{A}_{match}(\pi)-A(\pi) \stackrel{p}{\longrightarrow} 0.$$

(ii)  Under  Assumptions 1 and 2, then for any  $\pi \in \Pi$,  as $n\rightarrow\infty$, $$\left\{V^{\tau(X)}(\pi)+V^E(\pi)\right\}^{-1/2}\sqrt{n}\left\{\hat{A}_{match}(\pi)-B_M(\pi)-A(\pi)\right\} \stackrel{d}{\longrightarrow} N(0,1).$$
\end{theorem}

Theorem \ref{consistency and normality without bias correction} shows that the matching-based estimator $\amatch(\pi)$ converges in probability to the advantage function $A(\pi)$ regardless of the covariate dimension $p$.  Furthermore, $\amatch(\pi)$ is asymptotically normally distributed after subtracting the conditional bias term  $B_M(\pi)$. Importantly,  Theorem \ref{consistency and normality without bias correction} together with  Lemma \ref{order_conditional_bias}  lead to the following rate of convergence result.

\begin{corollary}\label{cor1}
 Under   Assumptions 1 and 2,  for any  $\pi \in \Pi$, 
$$\hat{A}_{match}(\pi)-A(\pi) = O_p\left(n^{-1/\max\{2, p\}}\right).$$
\end{corollary}
 
Therefore, when the covariate is scalar or bivariate,  $\hat{A}_{match}(\pi)$ is $\sqrt{n}$-consistent. However, the convergence rate can be very slow with a moderately large $p$, because $B_M(\pi)$ becomes the leading term. Therefore, it is imperative  to eliminate the conditional bias and thus remedy the slow rate for the matched-based estimator of the advantage function.

\subsection{Bias-Corrected Matching Estimator of the Advantage Function}\label{section Bias-corrected Matching-based Advantage Function}
We utilize the  idea of bias correction in \citet{abadie2011bias} to improve the matching estimator of the advantage function. Let $\lambda=\left(\lambda_1,\ldots,\lambda_p\right)^\top$ be a $p$-dimensional vector of nonnegative integers with the norm $\Vert \lambda \Vert_1= \sum_{j=1}^p \lambda_j$. Define that $x^\lambda= \prod_{j=1}^p (x_j)^{\lambda_j}$, where $x_j$ is the $j$th element of $p$-dimensional vector $x$. For the series $\{\lambda(j)\}_{j=1}^{\infty}$ containing all distinct vectors with nondecreasing $\Vert \lambda(j) \Vert_1$, define the $K$-dimensional vector $p^K(x)=\left(x^{\lambda(1)},\ldots,x^{\lambda(K)}\right)^\top$. The power series estimator \citep{newey1997series} of  $\mu(x,w)$ is given by
\begin{equation*}
\hat{\mu}(x,w)=\left(p^K(x)\right)^\top\left\{\underset{i:W_i=w}{\sum} p^K(X_i) \left(p^K(X_i)\right)^\top \right\}^{-} \underset{i:W_i=w}{\sum} p^K(X_i) Y_i,
\end{equation*}
where $(\cdot)^{-} $ is a generalized inverse. The bias-corrected matching estimator of the potential outcomes \citep{abadie2011bias} is 
\begin{equation*}
	\tilde{Y}_i(0)=\left\{ \begin{aligned}
		&Y_i, &\text{ if } W_i = 0,\\
		&\frac{1}{M} \underset{j \in \mathcal{J}_M(i)}{\sum} \left\{Y_j+\hat\mu(X_i,0)-\hat\mu(X_j,0)\right\}, &\text{ if } W_i = 1,
	\end{aligned} \right.
\end{equation*}
and
\begin{equation*}
	\tilde{Y}_i(1)=\left\{ \begin{aligned}
	&\frac{1}{M} \underset{j \in \mathcal{J}_M(i)}{\sum} \left\{Y_j+\hat\mu(X_i,1)-\hat\mu(X_j,1)\right\}, &\text{ if } W_i = 0,\\ &Y_i, &\text{ if } W_i = 1.
	\end{aligned} \right.
\end{equation*}
The bias-corrected  matching estimator of the ATE $\tau$ is  $\hat{\tau}_{bc}=n^{-1} \sum_{i=1}^n\{\Ytilde_i(1)-\Ytilde_i(1)\}$. Analogously, we  develop a bias-corrected matching-based advantage function:
\begin{equation}\label{Afun_bc_match}
	\hat{A}^{bc}_{match}(\pi)= \frac{1}{n}\sum^n_{i=1} \left(2\pi(X_i)-1\right) \left\{\tilde{Y}_i(1)-\tilde{Y}_i(0)\right\}.
\end{equation}
It follows  that 
$$\hat{A}^{bc}_{match}(\pi)=\hat{A}_{match}(\pi)-\hat{B}_M(\pi),$$
where the term
\begin{equation*}
	\hat{B}_M(\pi)=\frac{1}{n}\sum^n_{i=1} (2W_i-1)\left\{2\pi(X_i)-1\right\}\left[\frac{1}{M}\underset{j \in \mathcal{J}_M(i)}{\sum}\left\{\hat\mu(X_i,1-W_i)-\hat\mu(X_j,1-W_i)\right\}\right]
\end{equation*}
 serves as an estimator of $B_M(\pi)$.
The following assumption  contains regularity conditions for the power series estimator $\hat\mu(x,w)$ to ensure its fast convergence rate and thus eliminate the conditional bias. Detailed discussions can be found in  \citet{abadie2011bias} and \citet{otsu2017bootstrap}.

\begin{assumption}
Assume the following conditions:

(i) The support $\mathcal{X}_0$ is a Cartesian product of compact intervals.

(ii) $K=O(n^v)$, where $0 < v < \min\{2/(4p+3), 2/(4p^2-p)\}$.

(iii) For each $\lambda$ and $w = 0, 1$, the derivative $\partial^{\Vert \lambda \Vert_1} \mu(x,w)/\partial x^{\lambda_1}_1 \cdots \partial x^{\lambda_p}_p$ exists and its norm is upper bounded by $C^{\Vert \lambda \Vert_1}$ for some positive constant $C$.

\end{assumption}

We proceed to establish the  asymptotic normality of $\hat{A}^{bc}_{match}(\pi)$ in the following theorem.

\begin{theorem}\label{bias-corrected estimator normality}
	Under  Assumptions 1-3, then for any $\pi \in \Pi$, as $n\rightarrow\infty$, $$\left\{V^{\tau(X)}(\pi)+V^E(\pi)\right\}^{-1/2}\sqrt{n}\left\{\hat{A}^{bc}_{match}(\pi)-A(\pi)\right\} \stackrel{d}{\longrightarrow} N(0,1).$$
\end{theorem}
It follows immediately   that $\amatchbc(\pi) - A(\pi) = O_p(1/\sqrt{n})$, and thus the convergence rate of the bias-corrected matching estimator is faster when $p\geq 3$ than  that in Corollary \ref{cor1}.

\section{Policy Learning with Matching}\label{section implementation}
In the preceding discussion, we  introduced   matching estimators   of the advantage function $A(\pi)$ when the policy $\pi$ is fixed.  We now proceed to  policy learning with matching. We propose to learn the optimal policy, denoted as $\hat{\pi}$, by optimizing  the bias-corrected matching estimator of the advantage function  $\amatchbc(\pi)$ within a pre-determined  policy class $\Pi$:
    \begin{equation}
    \hat{\pi}=\arg \max_{\pi \in \Pi}\amatchbc(\pi).\label{learnedpolicy}
    \end{equation}

We  now proceed to describe the algorithm to solve the optimization problem \eqref{learnedpolicy}.  Furthermore, we derive a non-asymptotic regret bound for the learned policy $\hat{\pi}$.

\subsection{Implementation Detail}\label{implement}
Notice that $\amatchbc(\pi)$ has the following expression:
\begin{equation}\label{weighted-classification}
\amatchbc(\pi)=\frac{1}{n}\sum^n_{i=1} \left(2\pi(X_i)-1\right)\operatorname{sign}\left\{\tilde{Y}_i(1)-\tilde{Y}_i(0)\right\}\left|\tilde{Y}_i(1)-\tilde{Y}_i(0)\right|.
\end{equation}
Therefore, the problem of learning the optimal policy in \eqref{learnedpolicy} using the expression \eqref{weighted-classification} is equivalent to a weighted classification problem, where we  treat $\pi(\cdot)$ as the classifier to be trained, $\{\operatorname{sign}\{\tilde{Y}_i(1)-\tilde{Y}_i(0)\}\}^n_{i=1}$ as the labels, and $\{|\tilde{Y}_i(1)-\tilde{Y}_i(0)|\}^n_{i=1}$ as the weights.

Given a set of $n$ training samples $\{X_i,W_i,Y_i\}_{i=1}^n$ and a restricted policy class $\Pi$, we describe the implementation of MB-learning as follows:
\begin{enumerate}
	\item Impute the pairs of potential outcomes $\{Y_i(1),Y_i(0)\}^n_{i=1}$ using the bias-corrected matching estimators $\{\Ytilde_i(1),\Ytilde_i(0)\}^{n}_{i=1}$ defined in Section \ref{section Bias-corrected Matching-based Advantage Function}.
	\item  Construct the bias-corrected matching-based advantage function $\amatchbc(\pi)$ in  equation \eqref{weighted-classification}.
	\item Solve the optimization problem \eqref{learnedpolicy}   from the weighted classification perspective.
\end{enumerate}

The idea of recasting the optimization problem of the value function into a weighted classification framework has been previously discussed in the literature \citep[e.g.,][]{zhangb2012stat,zhao2012estimating,athey2021policy}. Directly solving this optimization problem in Step 3 is difficult since the estimated advantage function $\amatchbc(\pi)$ is not  convex. In this paper,  we employ the policy tree algorithm in    \citet{athey2021policy}, where we restrict the policy class $\Pi$ to be a class of fixed-depth decision trees and solve the weighted classification problem using  the R package \texttt{policytree} \citep{sverdrup2020policytree}. Due to its tree structure, the learned policy is easy to understand and interpret  for decision-makers. 

Alternatively, one may use some convex surrogate loss function to replace the nonsmooth indicator function in equation \eqref{weighted-classification}, thereby reformulating the policy learning problem into a convex optimization problem which can be solved  using off-the-shelf packages \citep{zhao2012estimating,liu2018augmented}. This is beyond the scope of this work.

\subsection{Regret Bound}\label{section Regret}
In the previous section, we  proved the asymptotic behaviour of the proposed matching-based estimators for the advantage function $A(\pi)$ for  a fixed policy $\pi$.   In this subsection, we  evaluate the performance of the learned  policy compared to the underlying optimal policy $\pi^*$ using the criterion of regret.  The utilitarian regret of a policy $\pi \in \Pi$ \citep{athey2021policy,zhouzhengyuan2023} is defined to be  $$R(\pi)= V(\pi^*) - V(\pi).$$ 
We proceed to investigate the convergence rate of the regret $R(\hat{\pi})$ with the learned matching-based policy $\hat{\pi}$ obtained from the optimization problem \eqref{learnedpolicy}. We aim to derive a non-asymptotic upper bound for $R(\hat{\pi})$  with high probability.  

For a given policy class $\Pi$, let $\text{VC}(\Pi)$ be its  Vapnik-Chervonenkis (VC) dimension \citep{vdV23}. In the following theorem, We assume that $\text{VC}(\Pi)$ is finite  to control the complexity of $\Pi$. Suppose that we learn the optimal policy by  MB-learning within a class of depth-$H$ decision trees, then the corresponding VC dimension is $\tilde{O}\left(2^H \log p\right)$ \citep{athey2021policy}. Here, the notation $\tilde{O}\left(f(n)\right)$  represents any factor upper bounded  by $g(f(n))f(n)$, where the function $g(\cdot)$ scales polylogarithmically with its argument. We also assume that $\Pi$ is countable, which is a standard technical condition to simplify the measurability  \citep[e.g.,][]{kitagawa2018Econometrica}. We write the notation $A \lesssim B$ to indicate that $A$ is less than $B$ up to a constant.

\begin{theorem}\label{non-asymptotic-bound}
	Under Assumptions 1-3, suppose that the outcome $|Y|\leq L$ almost surely for some positive constant $L$,  the policy class $\Pi$ is countable, and  $\text{VC}(\Pi)$ is finite, then for any $\delta>0$ and any $q > 1$, there exists an integer $N_{\delta,q}$ such that for all $n > N_{\delta,q}$, the regret $R(\hat{\pi})$ is upper bounded as follows with probability at least  $1-11\delta$:
	\begin{equation*}
	R(\hat{\pi}) \lesssim
	\frac{L}{\sqrt{M n^{1-1/q}}} \left(\sqrt{\text{VC}(\Pi)} + \sqrt{\frac{\log\frac{1}{\delta}}{M}}\right) +\frac{1}{\sqrt{n}}\left(\sqrt{\text{VC}(\Pi)}+\sqrt{\log\frac{1}{\delta}}\right)\left(\sqrt{V^*}  + L\right) + \frac{L \log\frac{1}{\delta}}{M n^{1-1/q}},
    \end{equation*}
where
$V^* =  4\sup_{\pi_1,\pi_2 \in \Pi} E\left[\left(\mu(X,1)-\mu(X,0)\right)^2 \left(\pi_1(X)-\pi_2(X)\right)^2 \right].$
\end{theorem}

It is equivalent to analyze the bound for $A(\pi^*)-A(\hat{\pi})$ in the proof. We decompose it into several terms, for which we derive  high-probability upper bounds, respectively. We remark that  it is challenging to derive the upper bound for the following term:
\begin{equation}\label{intractable term}
E\left[	\underset{\pi \in \Pi}{\sup}\left|\frac{1}{n}\sum^n_{i=1} (2W_i-1)\left\{\sum_{j=1}^n \pi(X_j)M_{ji}\right\}\epsilon_i\right| \right],
\end{equation}
because that matching induces pairwise correlations  between the items $\{\sum_{j=1}^n \pi(X_j)M_{ji}\}_{i=1}^n$. However, conditional on  $\{X_i,W_i,\epsilon_i\}_{i=1}^n$, it is   still eligible to utilize standard  empirical process techniques including symmetrization \citep{Boucheron2013ConcentrationInequalities} to show that the corresponding Rademacher process is sub-Gaussian. Consequently, we  leverage Dudley's entropy integral bound \citep{RomanVershynin2018HDP} and the order of $\max_{i \leq n}\KMi$ to derive an upper bound for the term \eqref{intractable term}.

Theorem \ref{non-asymptotic-bound} shows that the regret bound depends on the VC dimension $\text{VC}(\Pi)$, the number of matches $M$, and the quantity $V^*$. A better regret bound can be achieved if we decrease the complexity of the policy class or increase the number of matching. Similar discussion  can be found in \citet[Theorem 5]{abadie2006large}, where they stated that a larger $M$  generally leads to  a more efficient matching estimator of the ATE. Recalling that the unit-level term of $\bar{A}(\pi)$ is given by $\bar{A}_{i}(\pi)=\{2\pi(X_i)-1\}\{\mu(X_i,1)-\mu(X_i,0)\}$, the term  $V^*$ can be reformulated as $\sup_{\pi_1,\pi_2 \in \Pi} E\left[\bar{A}_{i}(\pi_1)-\bar{A}_{i}(\pi_2)\right]^2$. Therefore, we can interpret $V^*$ as the worst-case discrepancy between the advantage functions of any two policies in the class $\Pi$. See \citet{zhouzhengyuan2023} for related discussions.

 A direct consequence of Theorem \ref{non-asymptotic-bound}  is the convergence rate of the regret $R(\hat{\pi})$:
\begin{corollary}\label{cor2}
 Under  the assumptions in Theorem \ref{non-asymptotic-bound}, for any $q>1$, it holds that 
 $$R(\hat{\pi})=O_p\left(\sqrt{\frac{\text{VC}(\Pi)}{n^{1-1/q}}}\right).$$
\end{corollary}

\citet{athey2021policy} proved that the regret of AIPW-based policy learning with cross-fitting is rate-optimal with the order of $\sqrt{\text{VC}(\Pi)/n}$. Corollary \ref{cor2} suggests that the regret of the proposed MB-learning  has slower convergence rate for any $q>1$. Nevertheless, with  $q$  going to $\infty$,  the convergence rate  is almost of order  $\sqrt{\text{VC}(\Pi)/n}$.  That is,  MB-learning is almost rate-optimal.

Although  MB-learning  is not rate-optimal,  the simulation studies and real data analysis demonstrate the superiority of the proposed method compared to AIPW-based policy learning, especially when the sample size is not large or when the propensity score has extreme values.  

\citet{kitagawa2018Econometrica} showed that their IPW-based method has a regret bound faster than the order of $1/\sqrt{n}$ under the margin assumption, which is commonly used in classification literature \citep[e.g.,][]{tsybakov2004optimal,massart2006RiskBounds}. A faster convergence rate of the regret of MB-learning may be feasible under additional assumptions, but it is beyond the scope of this work.

\section{Simulation Studies}\label{section simulation}
\subsection{Data Generating Process}
We conduct extensive simulation studies to assess the finite sample performance of MB-learning. In each simulation setup, the covariate vector $X=(X_1,X_2,X_3,X_4)^\top$ is independently generated from $N(0,\mathbf{I}_4)$. The binary treatment $W$ is generated from the propensity score $\Pr(W=1|X)=\exp\left(l(X)\right)/\left[1+\exp\left(l(X)\right)\right]$, where we consider five different specifications of $l(X)$. Detailed information is presented in Table \ref{data_generating_table} of Supplementary Material. The outcome variable is generated by $Y = m(X) + Wc(X) + e, $ where $e$ is a noise term  generated from $N(0,1)$. We consider two scenarios of the main effect model $m(\cdot)$: (i). \textbf{Linear}: $m(X)=1 + 2 X_1 - X_2 + 0.5 X_3 - 1.5 X_4$; (ii). \textbf{Nonlinear}: $m(X)= 4\sin(X_1) + 2.5\cos(X_2) - X_3X_4$.   Moreover, We consider two scenarios of the contrast function $c(\cdot)$: (i). \textbf{Tree}: $c(X)=2I\left\{X_1>0, X_2>0\right\}-1$; (ii). \textbf{Non-tree}: $c(X)=2I\left\{2X_2-\exp(1+X_1)+2 >0\right\}-1$. We define the global optimal policy $\pi^{*}_{opt}$ as the  policy within all measurable functions from $\mathcal{X}$ to $\{0,1\}$ that maximizes the value function \citep{mo2021JASADistributionalRobustITR}, i.e., $\pi^{*}_{opt} = \arg\max_{\pi} V(\pi)$. Note that the contrast function $c(X)$ determines the structure of the global optimal policy since $\pi^{*}_{opt}(X)= I \left\{c(X)>0\right\}$. Therefore, when the contrast function has the tree structure,   the restricted policy class $\Pi$ includes the global optimal policy $\pi^{*}_{opt}$ in MB-learning. When  the contrast function is non-tree, the policy class $\Pi$ is misspecified, and  we investigate how the misspecification of $\Pi$ influences the proposed method. Similar simulation studies were conducted concerning tree-based  methods  in the literature \citep[e.g.,][]{zhangb2012stat,zhao2015treebasedITR,athey2021policy}.

\subsection{Compared Methods and Implementation Details}
The proposed method is compared with: (a) Q-learning \citep{qian2011performance}; (b) Augmented outcome weighted learning (AO-learning) \citep{liu2018augmented}, which improves O-learning with an argumentation term; (c) M-learning \citep{wu2020matched}, an alternative policy learning method based on matching; and (d)  AIPW-based policy learning with cross-fitting implemented by policy tree (denoted as Policy Tree) \citep{athey2021policy,sverdrup2020policytree}.

For MB-learning, we set the restricted policy class $\Pi$ as the class of depth-2 decision trees. Both  matching-based  and bias-corrected matching-based advantage functions in  \eqref{Afun_match} and  \eqref{Afun_bc_match} are considered, where the matching procedures are implemented by the R package \texttt{Matching}. We consider two different methods for bias-corrected MB-learning (MB-LR and MB-LASSO) to estimate $\mu(x,w)$. First, we  estimate $\mu(x,w)$ using linear regression  on $(1, X)$ based on units with $W=0, 1$ respectively (MB-LR) \citep{abadie2011bias}. The implementation details can be found in \citet[Page 299]{abadie2004Stata}.  We also employ  linear regression on $(1,X,X:X)$ with Lasso penalty  \citep{tibshirani1996regression} using the R package \texttt{glmnet}, where $X:X$ denotes all the squared and two-way interaction terms of $X$ (MB-LASSO). The number of matches $M$ is set to be $1 \text{ or } 5$, denoted by the suffix {-M1} or {-M5} respectively.  Therefore, we consider six variants of   MB-learning: MB-M1 and MB-M5, which utilize  the matching-based advantage function in  \eqref{Afun_match}; MB-LR-M1, MB-LR-M5, MB-LASSO-M1, and MB-LASSO-M5, which utilize  the  bias-corrected version in  \eqref{Afun_bc_match} with the outcome model $\mu(x,w)$  fitted by linear regression or Lasso.

For Policy Tree, we  derive the optimal policy within a class of depth-2 decision trees using the R package \texttt{policytree}, where the nuisance functions are estimated by the R function \texttt{double\_robust\_scores} in the R package \texttt{grf}.

For Q-learning, both parametric and non-parametric variants are considered. Specifically, we fit $Y$ on $(1,X,X:X,W, WX,W(X:X))$ with Lasso penalty using \texttt{glmnet}  for parametric Q-learning. For the non-parametric counterpart, we utilize random forest to fit $Y$ on $(X,W)$ by the R function \texttt{regression\_forest} in the R package {grf}. Here, we leverage the suffixes {-Para} and {-RF} to denote the parametric and nonparametric variants respectively.

For AO-learning, the parametric and non-parametric variants are also considered, represented by the suffixes {-Para} and {-RF} respectively.  We estimate the propensity score using logistic regression for $W$ on $(1,X)$ using the R function \texttt{glm} in the parametric variant, while \texttt{regression\_forest} is used for the non-parametric counterpart. The optimal policy is obtained using the R package \texttt{WeightSVM} with Gaussian kernel where the tuning parameters are chosen via five-fold cross-validation. In addition, for the augmentation term defined in \citet{liu2018augmented}, we fit $Y$ on $X$ employing \texttt{regression\_forest}.

For M-learning, we construct the matched sets defined in \citet{wu2020matched} using \texttt{Matching} and learn the optimal policy using \texttt{WeightSVM}, where the kernel and tuning parameters are chosen using the same configuration as AO-learning.  We only consider 1-NN matching (NN-matching with $M=1$), because  our preliminary simulation studies and  the numerical results in \citet{wu2020matched} indicate that the simplest 1-NN matching  yields  better performance with less computational burden compared to $M$-NN matching with $M>1$.

For each scenario, we vary the size of training dataset $n \in \{200,500,1000\}$. We use the empirical value function $\mathbb{P}_n[Y(\pi(X))]$ as the  performance criterion, where $\mathbb{P}_n$ is the empirical average over an independent testing dataset with size $20000$. Each simulation setup is replicated 200 times.

\subsection{Simulation Results}

\begin{figure}[t]
	\centering
	\includegraphics[width=1\linewidth]{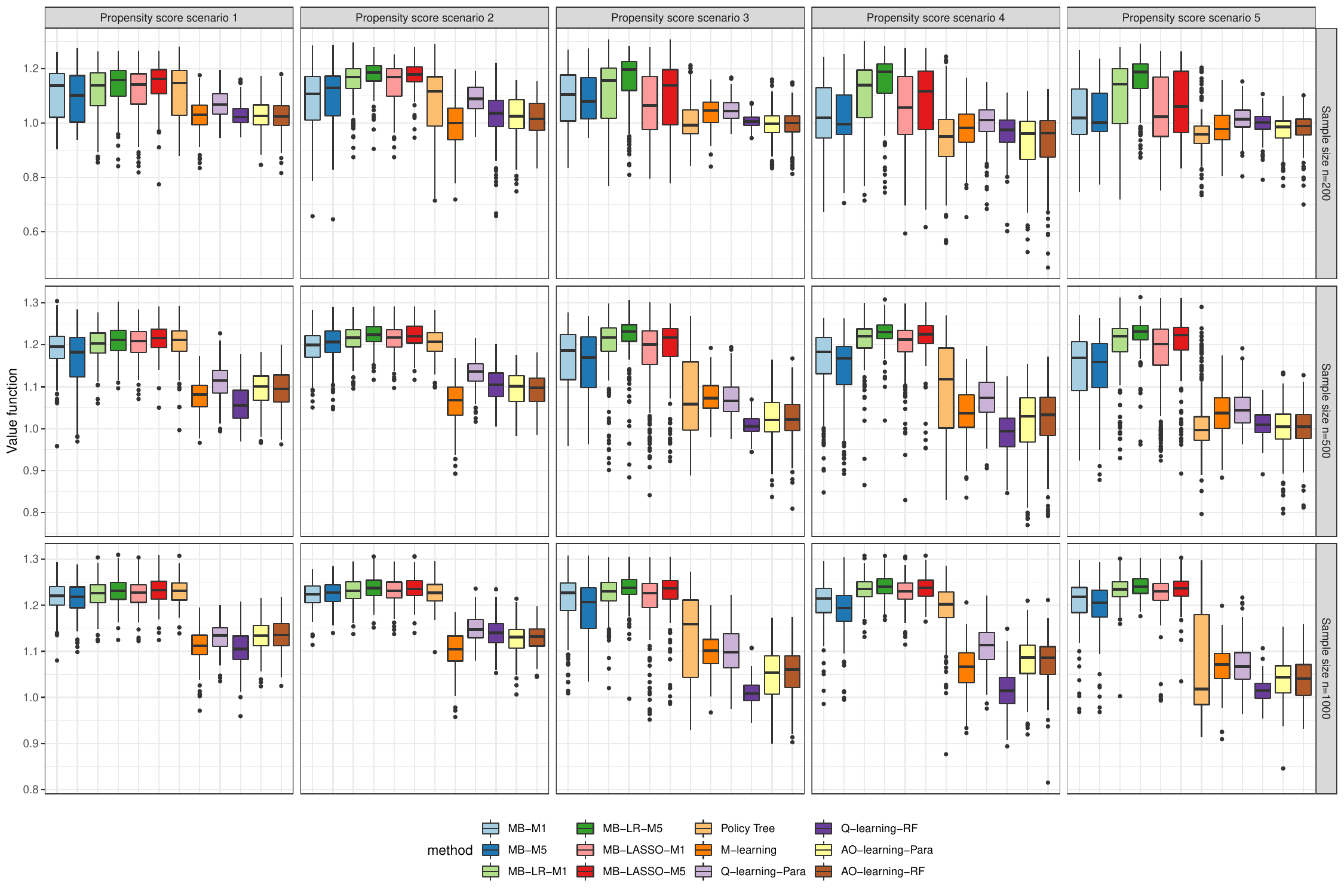}
	\caption{Boxplot of empirical value functions  where the main effect model is \textbf{linear} and the contrast function is \textbf{tree}. The global optimal value function is 1.25. The three rows represent  different sample sizes $n=200, 500$, or 1000. The five columns represent  different propensity score models  presented in Table \ref{data_generating_table} of Supplementary Material. }
	\label{performance-linearmain-tree}
\end{figure}

We present the simulation results in this subsection. Figure \ref{performance-linearmain-tree}  illustrates the performance in five propensity score scenarios with three different sample sizes, where the main effect model is linear and the contrast function is tree. The global optimal value function $V(\pi^{*}_{opt})$ is 1.25 in this setup. All the approaches exhibit more satisfactory performance with a larger sample size.  As expected, the six variants of MB-learning are competitive since the policy $\pi^{*}_{opt}$ is included in the class $\Pi$. Among them, MB-LR-M5 and MB-LASSO-M5 correctly fit the main effect model and   yield the  best performance, demonstrating the necessity of bias correction. The results indicate that MB-learning with a larger number of matches $M$ performs slightly better. However, preliminary simulations suggest that MB-learning is insensitive to $M>5$. Notably, Policy Tree yields the best performance among the remaining methods. Nevertheless, although Policy Tree  behaves comparably as MB-learning  in Scenarios 1 and 2, its performance deteriorates in Scenarios 3-5 with significant variability due to that the two treatment groups are extremely imbalanced. These phenomena align with  discussions in the literature on AIPW-based approaches  \citep[e.g.,][]{kangschafer2007,zhao2019efficient}. In contrast, the variants of MB-learning mitigate this issue by circumventing the use of propensity score in estimation and learning and  exhibit appealing robustness.

\begin{figure}[t]
	\centering
	\includegraphics[width=1\linewidth]{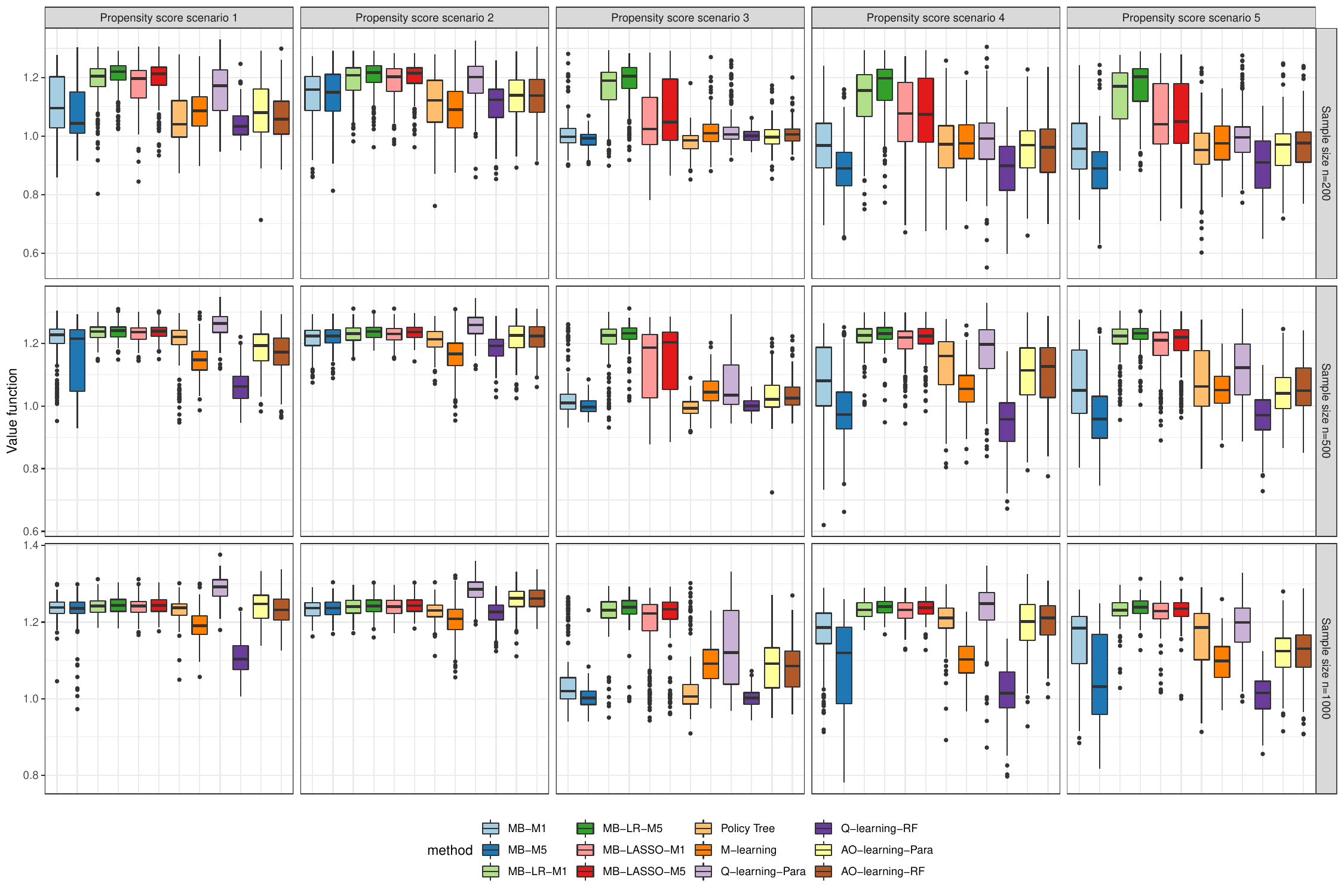}
	\caption{Boxplot of empirical value function  where the main effect model is \textbf{linear} and the contrast function is \textbf{non-tree}. The global optimal value function  is 1.36. The three rows represent  different sample sizes $n=200, 500$, or 1000. The five columns represent  different propensity score models  presented in Table \ref{data_generating_table} of Supplementary Material.}
	\label{performance-linearmain-nottree}
\end{figure}

Figure \ref{performance-linearmain-nottree} presents the results with linear main effect model and non-tree contrast function, where the global optimal value function  is 1.36. These scenarios explore the impact of misspecifying the  policy class. In general, MB-learning with bias correction and Q-learning-Para outperform the remaining methods. Among them, Q-learning-Para dominates the others in Scenarios 1 and 2,  which can be explained by the reasonable accommodation of the outcome mean model $\mu(x,w)$. However, in Scenarios 3 and 5, the variants of MB-learning with bias correction exhibit   superior performance, while Q-learning-Para deteriorates in these situations due to the extreme proportion of the two treatment groups particularly when  the sample size is small ($n=200$). These results demonstrate the competitive performance of MB-learning, which can be attributed to the desirable approximation of fixed-depth decision tree to the global optimal policy $\pi^{*}_{opt}$ and the robustness of matching. Among (A)IPW-based approaches, Policy Tree and AO-learning have favorable performance in Scenarios 1 and 2, while they behave unsatisfactorily in the remaining settings. It is noteworthy that M-learning yields similarly unfavorable performance as  MB-learning without bias correction, which may  be partly explained by the matching bias.

Figures \ref{performance-nonlinearmain-tree} and \ref{performance-nonlinearmain-nottree} in Supplementary Material display the results with nonlinear main effect model where the contrast function is tree or non-tree. Generally, MB-LASSO has the best performance. In contrast, MB-LR behaves similarly to MB-LASSO in Scenarios 1, 2, and 4, whereas the performance  is relatively poor in Scenarios 3 and 5. This is not surprising because that the conditional  bias is relatively minor in the former scenarios, as illustrated by the results of MB-learning without bias correction.  The conditional  bias becomes more significant in the latter situations. Hence, bias correction based on linear regression may be insufficient to eliminate   bias in Scenarios 3 and 5. In addition, Q-learning-RF is competitive in these settings, outperforming its parametric counterpart as is expected.  The remaining approaches exhibit similar performance  as those in Figures \ref{performance-linearmain-tree} and \ref{performance-linearmain-nottree}. 

We conduct  more numerical studies  comparing MB-learning with Policy Tree when $n$ varies from 200 to 8000. The results are displayed in Figures \ref{compare-with-policytree-trueboundary-tree} and \ref{compare-with-policytree-trueboundary-nottree} in Supplementary Material. MB-learning with bias correction outperforms Policy Tree by a large margin when the sample size is small, whereas the performance becomes nearly the same when the sample size increases to 8000. This empirical evidence agrees with the theoretical result of convergence rate in Corollary \ref{cor2}. This phenomenon can be attributed to the cross-fitting technique employed in Policy Tree: the nuisance functions are fitted with random forests, which may be inadequate with a small sample size.  Similar observations  were made in \citet{wu2022matching}.  Finally, Policy Tree encounters the ``phase transition'' phenomenon in Scenarios 3 and 5, whereas MB-learning does not have this issue in the simulations.

\section{Application: Treatment Allocation in the NSW Training Program}\label{section application}
The  NSW Program aimed to provide work experience for job-seekers encountering  economic and social difficulties before enrollment. The individuals were randomly assigned to the  training program-exposed group or the unexposed group. The causal effect of the training program on post-treatment earnings has been analyzed in many literature \citep[e.g.,][]{lalonde1986,dehejia1999causal,dehejia2002propensity,abadie2011bias,imai2014covariate}. We  utilize the NSW datasets, available at \href{https://users.nber.org/~rdehejia/nswdata2.html}{\textcolor{black}{https://users.nber.org/~rdehejia/nswdata2.html}}, to demonstrate the performance of MB-learning in identifying the optimal policy.

\begin{table}[t]
	\centering
	\caption{The normalized differences of covariates in different datasets.}
	\label{normalizeddifference_table}
	\resizebox{0.7\textwidth}{!}{
		\begin{tabular}{ccccccccc}
			\toprule
			& \multicolumn{8}{c}{Covariate} \\
			\cmidrule{2-9}
			Dataset & Age & Education & Black & Hispanic & Married & No degree & RE74 & RE75 \\
			\midrule
			DW & 0.107 & 0.141 & 0.044 & -0.175 & 0.094 & -0.304 & -0.002 & 0.084 \\
			DW-CPS3 & -0.128  &  0.074 &  0.930  & -0.240 & -0.426 &   0.047 &  -0.378 &  -0.148   \\
			\bottomrule
	\end{tabular}}
\end{table}

We first implement the proposed method using the dataset used in \citet{dehejia1999causal} (DW), which is a randomized experiment including 185 exposed subjects ($W=1$) and 260 unexposed subjects ($W=0$). The eight-dimensional pre-treatment covariates $X$ are considered as confounders, which consist of age, education,   race (1 if black, 0 otherwise), ethnicity  (1 if Hispanic, 0 otherwise), marital status (1 if married, 0 otherwise), academic degree (1 if no degree, 0 otherwise), earnings in 1974 (RE74), and earnings in 1975 (RE75). The outcome of interest $Y$ is the earnings in 1978 (RE78). To assess  covariate balance between the treatment groups,  we calculate the normalized  difference \citep{abadie2011bias}, defined as $(\bar{X}_1-\bar{X}_0)/\sqrt{(S^2_1+S^2_0)/2}$, where $\bar{X}_w = \sum_{i:W_i=w} X_i/n_w$, $S^2_w=  \sum_{i:W_i=w}(X_i-\bar{X}_w)^2/(n_w-1)$, and $n_w$ is the number of units in group $w$, $w=0,1$. The  normalized  differences  are displayed in Table \ref{normalizeddifference_table}, indicating  no significant difference in covariate distributions between the groups.

We proceed to derive the optimal policy $\hat{\pi}$ that assigns personalized treatments to the job-seekers with the goal of maximizing the post-treatment earnings. Due to ethical constraints, the black and Hispanic variables   are  not considered in learning the optimal policy, although they are used in estimating the advantage function or the value function.

We apply the six variants of MB-learning  discussed in Section \ref{section simulation} to the DW dataset and focus on the comparison with Policy Tree in \citet{athey2021policy}, where the optimal policies by MB-learning and Policy Tree are both estimated within depth-2 decision trees using the  R package \texttt{policytree}. The estimated policy $\hat{\pi}$ is evaluated by an AIPW estimator of the value function $V(\hat{\pi})$: $$\mathbb{P}_n\left[\frac{YI\left\{W=\hat{\pi}(X)\right\}}{\tilde{e}(\hat{\pi}(X),X)} - \frac{I\left\{W=\hat{\pi}(X)\right\}-\tilde{e}(\hat{\pi}(X),X)}{\tilde{e}(\hat{\pi}(X),X)}  \tilde{\mu}(X,\hat{\pi}(X))\right],$$ where $\tilde{e}(w,x)$ and $\tilde{\mu}(x,w)$ represent the estimators of $e(w,x)=\Pr(W=w|X=x)$ and $\mu(x,w)$ respectively. Notably, since this dataset is collected from a randomized experiment, $\tilde{e}(w,x)$ is attained by the empirical proportion of the exposed or unexposed group. The  $\tilde{\mu}(x,w)$ is derived by random forest. Analogous to the arguments in \citet{zhao2015doubly}, we conduct a cross-validated analysis to assess the performance of the estimated policy. Specifically, the dataset is randomly split into five parts with equal sample sizes, where four parts are used to estimate the optimal policy and the remaining  part is used to  evaluate the performance. The training and testing sets are permuted five times such that each part is used as the testing dataset in turn. Cross-validated values for each approach are obtained by averaging the results across the five testing sets. This procedure is repeated  100 times. To  evaluate the effectiveness of individualized policy assignments, we also explore the one-size-fits-all policy that all units are either exposed or unexposed, where we calculate the AIPW-based adjusted mean value.

Table \ref{crossvalidated_table} reports the cross-validated results of the value functions. The results indicate that the universal policy of always being exposed is  more beneficial  for  job-seekers compared to MB-learning and Policy Tree. One possible explanation is that the job training program generally increases earnings for the  population, which implies limited heterogeneity. This phenomenon was discussed in \citet{zhouzhengyuan2023}. The results may be partly attributed to   misspecification of the restricted policy class, which  lead to  deteriorated  performance of MB-learning and Policy Tree as discussed in Section \ref{section simulation}.  Among the personalized policies,  MB-LASSO-M5 exhibits the best performance, while the other variants of MB-learning are competitive. In contrast, Policy Tree performs the worst in terms of the cross-validated values, which can be explained by  the usage of machine learning in estimating nuisance functions with a relatively small sample. The tree structures of the estimated optimal policies  by MB-LASSO-M5 and Policy Tree are drawn in Figure \ref{learned-policy-DW} in Supplementary Material.

Moreover, we investigate the performance of MB-learning on observational data. The DW dataset is composited with additional 429 unexposed units from Current Population Survey (CPS3), resulting in a synthetic observational study (DW-CPS3).  The normalized  difference values of the DW-CPS3 dataset presented in Table \ref{normalizeddifference_table} suggest that the covariate discrepancy between the exposed and unexposed groups is not ignorable particularly for race, marital status, and earnings. We conduct a cross-validated analysis in this observational study, where both $\tilde{e}(w,x)$ and $\tilde{\mu}(x,w)$ are estimated using random forest. As shown in Table \ref{crossvalidated_table}, the variants of MB-learning have superior performance. Among them,  MB-M5 yields the highest cross-validated value and the smallest standard deviation. It is noteworthy that MB-learning slightly deteriorates after bias correction. One possible reason is the unsatisfactory fit of $\mu(x,w)$, which may be highly nonlinear in this population. In addition, the results illustrate that the individualized policies learned by MB-learning outperform the universal policy of always being exposed. These results are likely to be related to the presence of heterogeneity in this composite dataset, given that the CPS3 samples are nationally representative when they were collected across the United States \citep{smith2005does}. The job training program may not universally benefit these composite samples, as they were drawn from different local labor markets. On the other hand, Policy Tree is dominated by the variants of MB-learning, exhibiting  the highest variability and yielding a cross-validated value  worse than the universal policy of always being exposed. We  present the tree-structured optimal policies learned by MB-M5 and Policy Tree in Figure \ref{learned-policy-DW-CPS3} in Supplementary Material.

\begin{table}[t]
	\centering
	\caption{Mean cross-validated values in different datasets. The standard deviation is in the parenthesis. In DW dataset, the AIPW estimators of the value function of  the ``always being exposed" policy and ``always being unexposed" policy are 6244.9 and 4617.5, respectively. In DW-CPS3 dataset, the AIPW estimators of the value function of  the ``always being exposed" policy and ``always being unexposed" policy are 6614.2 and 5781.6, respectively.}
	\label{crossvalidated_table}
	\resizebox{\textwidth}{!}{\begin{tabular}{cccccccc}
		\toprule
		& \multicolumn{7}{c}{Policy learning method} \\
		\cmidrule{2-8}
		Dataset & MB-M1 & MB-M5 & MB-LR-M1 & MB-LR-M5 & MB-LASSO-M1 & MB-LASSO-M5 & Policy Tree \\
		\midrule
		DW & 5936.7 (220.5) & 6036.4 (199.8) & 5898.3 (263.8) & 5956.0 (198.8) & 5916.7 (214.8) & \textbf{6038.7} (226.2) & 5839.8 (226.4)  \\
		DW-CPS3 & 6725.4 (167.5)  & \textbf{6959.0} (127.8) &  6757.1 (154.0)  & 6909.7 (163.4) & 6719.0 (159.8) &   6949.8 (135.8) & 6594.6 (206.0)  \\
		\bottomrule
	\end{tabular}}
\end{table}

\section{Discussion}\label{section discussion}
In this paper, we proposed a framework of matching-based policy learning for identifying the optimal policy in observational studies by  leveraging various matching estimators in the  causal inference literature. Specifically, we developed a (bias-corrected) matching-based estimator of the advantage function and estimated the optimal policy by directly maximizing the resulting estimate over a class of policies. The proposed method possesses  robustness in observational studies by circumventing the use of inverse probabilities. Meanwhile, covariate balancing is naturally attained by matching. Theoretically, we established a non-asymptotic bound for the regret of the learned policy, providing rigorous performance guarantee for MB-learning. Empirically, we demonstrated the competitive performance of the proposed method, particularly when  extreme weights are present or when the sample size is not large.

There are several interesting directions in which this work can be further developed. First, MB-learning is intended for low-dimensional settings with binary treatment. In causal inference literature, \citet{antonelli2018doubly} proposed a doubly robust matching estimator to handle high-dimensional covariates, where variable selection techniques were used for matching with estimated propensity score and prognostic score. In addition, \citet{wu2022matching} developed a matching estimator of average causal exposure-response function based on generalized propensity score. These estimators offer inspiration for adapting MB-learning to different scenarios. Second,  we developed MB-learning with a fixed number of matches $M$, and  our discussion in Section \ref{section Regret} demonstrated that a larger $M$ yields a lower regret bound. \citet{lin2023estimation} established an insightful result that the bias-corrected matching estimator \citep{abadie2011bias} with a diverging $M$ exhibits doubly robustness and semiparametric efficiency under regularity conditions.  MB-learning with a diverging $M$ may yield a doubly robust policy learning approach with an improved regret bound. The robustness of MB-learning may be partly attributed to  covariate balance between the treatment groups, which is achieved through directly matching on covariates. It  would be of considerable
interest to integrate more sophisticated covariate balancing techniques for learning the optimal policy \citep[e.g.,][]{wong2018kernel,wang2020minimal,fan2022optimal,DaiYan2024}. Finally, extension of MB-learning to deal with multicategory or continuous policy is of future research interest \citep[e.g.,][]{LiXuqiao2025}.

	\bibliographystyle{ecta}
	\begin{spacing}{1.2}
		\bibliography{paper-ref}
	\end{spacing}

\appendix
\newpage

\setcounter{figure}{0}
\setcounter{table}{0}
\renewcommand{\thefigure}{S\arabic{figure}}
\renewcommand{\thetable}{S\arabic{table}}

\noindent \textbf{\Large  Supplementary Material}

\section{Technical Proofs}
Before the proofs, we introduce some notations here. Let $(\Theta,d)$ be a semi-metric space, that is, the space $\Theta$ is equipped with a semi-metric $d$. For any $\epsilon>0$, let  $N(\epsilon,\Theta,d)$ denote the covering number of  $\Theta$, i.e., the minimal number of $d$-balls with radius $\epsilon$ needed to cover $\Theta$. We call a set of points is $\epsilon$-separated if the distance between every pair of points is strictly larger than $\epsilon$, and let $D(\epsilon,\Theta,d)$ denote the packing number of $\Theta$, i.e., the maximum number of points in $\Theta$ such that they are $\epsilon$-separated. The entropy numbers are the logarithms of covering and packing numbers respectively. These definitions can be found elsewhere \citep[e.g.,][]{vdV23}.
Throughout the proofs, we use $C$ to denote a generic positive constant, which may change from place to place.
\subsection{Proof of Lemma 1}
\begin{proof}
Let $\mathbf{Y}=(Y_1,\ldots,Y_n)^\top$. Note that $$\hat{Y}(1)-\hat{Y}(0)=(2W_i-1)\left(Y_i- \frac{1}{M} \underset{j \in \mathcal{J}_M(i)}{\sum} Y_j\right)=(2W_i-1)\left(Y_i - \frac{1}{M} M^\top_i \mathbf{Y} \right),$$  and that
\begin{equation*}
\begin{aligned}
\sum_{i=1}^n \left(2\pi(X_i)-1\right)(2W_i-1) M^\top_i \mathbf{Y}=&\sum_{i=1}^n  \left(2\pi(X_i)-1\right)(2W_i-1) \sum_{l=1}^n M_{il} Y_l \\ =&  \sum_{l=1}^n \left[\sum_{i=1}^n \left(2\pi(X_i)-1\right)(2W_i-1)  M_{il} \right]Y_l \\ =& \sum_{l=1}^n \left[\underset{i: M_{il}=1}{\sum} \left(2\pi(X_i)-1\right)(1-2W_l) \right] Y_l \\ =& -\sum_{l=1}^n (2W_l-1) K_M(\pi,l) Y_l,
\end{aligned}
\end{equation*}
where the third equality follows from the fact that, if $M_{il}=1$, then $l \in \mathcal{J}_M(i)$ and thus  $W_i = 1-W_l$ by the definition of $\JMi$.
Combining these results, we have
\begin{equation*}
\begin{aligned}
	\amatch(\pi)=& \nmeani \left(2\pi(X_i)-1\right)\left(\Yhat_i(1)-\Yhat_i(0)\right) \\
	=& \nmeani \left(2\pi(X_i)-1\right)(2W_i-1)\left(Y_i - \frac{1}{M} M^\top_i \mathbf{Y} \right)  \\
	=& \nmeani (2W_i-1) \left[\left(2\pi(X_i)-1\right) + \frac{K_M(\pi,i)}{M}\right] Y_i.
\end{aligned}
\end{equation*}
\end{proof}

\subsection{Proof of Lemma 2}
\begin{proof}
Noting that $| \KMpi | \leq  \sum_{j=1}^n M_{ji} = \KMi$,  the lemma follows directly from the Lemma 3 in \citet {abadie2006large}.
\end{proof}

\subsection{Proof of Lemma 3}
\begin{proof}
Using the same argument in \cite{abadie2002simple}, we define two $n \times n$ matrices $\mathbf{G}_1$ and  $\mathbf{G}_0$:
\begin{equation*}
	\mathbf{G}_{1,ij}=\begin{cases}
		1& \text{ if } W_i = 1 \text{ and } i=j,\\
		\frac{1}{M}& \text{ if } W_i = 0 \text{ and } j \in \JMi,\\
		0& \text{ otherwise},
	\end{cases}
    \text{ and }
	\mathbf{G}_{0,ij}=\begin{cases}
		1& \text{ if } W_i = 0 \text{ and } i=j,\\
		\frac{1}{M}& \text{ if } W_i = 1 \text{ and } j \in \JMi,\\
		0& \text{ otherwise}.
	\end{cases}
\end{equation*}
Let $\mathbf{\hat{Y}}(1)=\left(\Yhat_1(1),\ldots,\Yhat_n(1)\right)^\top$ and $\mathbf{\hat{Y}}(0)=\left(\Yhat_1(0),\ldots,\Yhat_n(0)\right)^\top$ denote the $n \times 1$ vectors of matching estimators for potential outcomes under each treatment. Then, $$\mathbf{\hat{Y}}(1)=\mathbf{G}_1 \mathbf{Y} \text{ and } \mathbf{\hat{Y}}(0)=\mathbf{G}_0 \mathbf{Y}.$$
Furthermore, let $\bm{\mu}\mathbf{(X \text{,} W)}=\left(\mu(X_1,W_1),\ldots,\mu(X_n,W_n)\right)^\top$. It follows that
\begin{equation}\label{eq_Gmu}
	\mathbf{G_1}\bm{ \mu}\mathbf{(X \text{,} W)}= \mathbf{G_1}\bm{ \mu}\mathbf{(X \text{,} 1)} \text{ and } \mathbf{G_0}\bm{ \mu}\mathbf{(X \text{,} W)}= \mathbf{G_0}\bm{ \mu}\mathbf{(X \text{,} 0)}
\end{equation} by the definition of $\mathbf{G}_1$ and $\mathbf{G}_0$.

Let $\bm{\phi}$ be the $n \times 1$ vector taking value in $\{-1,1\}^n$ with the $i$th element  $2\pi(X_i)-1$,  and $\bm{\epsilon}$ be the  $n \times 1$ vector with the $i$th element $\epsilon_i$ which is equal to $Y_i-\mu(X_i,W_i)$. Then, we obtain
\begin{equation}\label{eq_decomposition}
\begin{aligned}
	\amatch(\pi)=&  \frac{1}{n} \bm{\phi}^\top (\mathbf{G_1}-\mathbf{G_0}) \mathbf{Y} \\
	=& \frac{1}{n} \bm{\phi}^\top  (\mathbf{G_1}-\mathbf{G_0}) \bm{\mu}\mathbf{(X \text{,} W)} + \frac{1}{n} \bm{\phi}^\top  (\mathbf{G_1}-\mathbf{G_0}) \bm{\epsilon} \\
	=& \frac{1}{n} \bm{\phi}^\top  \mathbf{G_1} \bm{\mu}\mathbf{(X \text{,} 1)} - \frac{1}{n} \bm{\phi}^\top  \mathbf{G_0} \bm{\mu}\mathbf{(X \text{,} 0)} + \frac{1}{n} \bm{\phi}^\top  (\mathbf{G_1}-\mathbf{G_0}) \bm{\epsilon} \\
	=& \frac{1}{n} \bm{\phi}^\top \left[\muXOne - \muXZero\right] + \frac{1}{n} \bm{\phi}^\top  (\mathbf{G_1}-\mathbf{G_0}) \bm{\epsilon} \\
	 & +  \frac{1}{n} \bm{\phi}^\top  (\mathbf{G_1}-\mathbf{I_n}) \bm{\mu}\mathbf{(X \text{,} 1)} - \frac{1}{n} \bm{\phi}^\top  (\mathbf{G_0}-\mathbf{I_n}) \bm{\mu}\mathbf{(X \text{,} 0)},
\end{aligned}
\end{equation}
where the third equality follows from equation \eqref{eq_Gmu} and that $\mathbf{I_n}$ is the $n \times n$ identity matrix.
Note that the  elements of the two $n \times n$ matrices $\mathbf{G_1}-\mathbf{I_n}$ and  $\mathbf{G_0}-\mathbf{I_n}$ are given by
\begin{equation*}
	(\mathbf{G_1}-\mathbf{I_n})_{ij}=\begin{cases}
		-1& \text{ if } W_i = 0 \text{ and } i=j,\\
		\frac{1}{M}& \text{ if } W_i = 0 \text{ and } j \in \JMi,\\
		0& \text{ otherwise},
	\end{cases}
     \text{ and }
	(\mathbf{G_0}-\mathbf{I_n})_{ij}=\begin{cases}
		-1& \text{ if } W_i = 1 \text{ and } i=j,\\
		\frac{1}{M}& \text{ if } W_i = 1 \text{ and } j \in \JMi,\\
		0& \text{ otherwise}.
	\end{cases}
\end{equation*}
For each row $i$ such that $W_i = 0$, we have $(\mathbf{G_1}-\mathbf{I_n})_{i \cdot}\muXOne=\frac{1}{M} \sum_{j \in \JMi} \mu(X_j,1) - \mu(X_i,1)$, where $(\mathbf{G_1}-\mathbf{I_n})_{i \cdot}$ is the $i$th row of  $(\mathbf{G_1}-\mathbf{I_n})$; for each row $i$ such that $W_i = 1$, we have  $(\mathbf{G_1}-\mathbf{I_n})_{i \cdot}\muXOne=0$. Hence,
\begin{equation*}
	\bm{\phi}^\top  (\mathbf{G_1}-\mathbf{I_n}) \bm{\mu}\mathbf{(X \text{,} 1)} = \underset{i: W_i=0}{\sum} \left(2\pi(X_i)-1\right) \left[\frac{1}{M} \underset{ j \in \JMi}{\sum} \mu(X_j,1) - \mu(X_i,1)\right],
\end{equation*}
and
\begin{equation*}
	\bm{\phi}^\top  (\mathbf{G_0}-\mathbf{I_n}) \bm{\mu}\mathbf{(X \text{,} 0)} = \underset{i: W_i=1}{\sum} \left(2\pi(X_i)-1\right) \left[\frac{1}{M} \underset{ j \in \JMi}{\sum} \mu(X_j,0) - \mu(X_i,0)\right]
\end{equation*}
by  similar arguments. Furthermore, it follows that
\begin{equation}\label{eq_decomposition_term3}
	\begin{aligned}
			\bm{\phi}^\top  (\mathbf{G_1}-\mathbf{I_n}) \bm{\mu}\mathbf{(X \text{,} 1)} - &\bm{\phi}^\top  (\mathbf{G_0}-\mathbf{I_n}) \bm{\mu}\mathbf{(X \text{,} 0)} \\
			= &\sum^n_{i=1} (2W_i-1)\left(2\pi(X_i)-1\right)\frac{1}{M}\underset{j \in \mathcal{J}_M(i)}{\sum}\left(\mu(X_i,1-W_i)-\mu(X_j,1-W_i)\right)
	\end{aligned}
\end{equation}
After some algebra and  similar arguments as the proof of Lemma \ref{linear_combination}, we have
\begin{equation}\label{eq_decomposition_term1}
	\bm{\phi}^\top \left[\muXOne - \muXZero\right] = \sum^n_{i=1}\left(2\pi(X_i)-1\right) \left(\mu(X_i,1)-\mu(X_i,0)\right),
\end{equation}
and
\begin{equation}\label{eq_decomposition_term2}
 \bm{\phi}^\top  (\mathbf{G_1}-\mathbf{G_0}) \bm{\epsilon}
	= \sum^n_{i=1} (2W_i-1)\left[(2\pi(X_i)-1)+\frac{K_M(\pi,i)}{M} \right]\epsilon_i.
\end{equation}
Combining \eqref{eq_decomposition}, \eqref{eq_decomposition_term3}, \eqref{eq_decomposition_term1} and \eqref{eq_decomposition_term2}, we complete the proof.
\end{proof}

\subsection{Proof of Lemma 4}
\begin{proof}
To prove the lemma, we need to utilize the order of the distances between unit $i$ and its matches. Let $j_m(i), m=1,\ldots,M$ denote the $m$th match to unit $i$, then we have $\JMi=\left\{j_1(i),\ldots,j_M(i)\right\}$. It follows that  $W_{j_m(i)}=1-W_i$ and that
\begin{equation*}
	\underset{k: W_k=1-W_i}{\sum} I\left\{ \Vert X_k-X_i \Vert\leq \Vert X_{j_m(i)}-X_i\Vert \right\} = m,
\end{equation*}
where, as in \citet{abadie2006large}, we ignore the ties in matching.  We reexpress $B_M(\pi)$ in \eqref{eq_BMpi} using the indices $\{j_m(i), m=1,\ldots,M; i=1,\ldots,n\}$ as follows:
\begin{equation}\label{BMpi_rewrite}
	\begin{aligned}
	B_M(\pi)=&\frac{1}{n}\sum^n_{i=1} (2W_i-1)\left(2\pi(X_i)-1\right)\frac{1}{M} \sum_{m=1}^M \left[\mu(X_i,1-W_i)-\mu(X_{j_m(i)},1-W_i)\right] \\
	=& \frac{1}{nM} \sum_{i=1}^n \sum_{m=1}^M \left(2\pi(X_i)-1\right)\left\{W_i\left[\mu(X_i,0)-\mu(X_{j_m(i)},0)\right] - (1-W_i)\left[\mu(X_i,1)-\mu(X_{j_m(i)},1)\right]\right\} \\
	\equiv& \frac{1}{nM} \sum_{i=1}^n \sum_{m=1}^M B^{i,m}_{M}(\pi).
	\end{aligned}
\end{equation}
Under the assumption that $\mu(x,0)$ and  $\mu(x,1)$ are Lipschitz, we have $\left| B^{i,m}_{M}(\pi) \right| \leq C \Vert U_{i,m} \Vert$ for some positive constant $C$, where $U_{i,m}=X_i-X_{j_m(i)}$ is the matching discrepancy between unit $i$ and $j_m(i)$. It follows from Cauchy-Schwarz inequality that
\begin{equation*}
	E\left[n^{2/p} \left(B_M(\pi)\right)^2\right] \leq n^{2/p} E \left[ \frac{1}{nM} \sum_{i=1}^n \sum_{m=1}^M \left(B^{i,m}_{M}(\pi)\right)^2 \right] \leq  C^2 n^{2/p} E \left[ \nmeani \Vert U_{i,M} \Vert^2\right].
\end{equation*}

Using the same argument in \citet[Proof of Theorem 1, Page 258]{abadie2006large}, we have $C^2 n^{2/p} E \left[ \nmeani \Vert U_{i,M} \Vert^2\right]=O(1)$. Therefore,  it holds that $ B_M(\pi)=O_p\left(n^{-1/p}\right)$ by Markov's inequality.

\end{proof}

\subsection{Proof of Theorem 1}
\subsubsection{Proof of Theorem 1 (i)}
\begin{proof}
By Lemma \ref{lemma_Amatch_decomposition}, we rewrite $$\hat{A}_{match}(\pi) - A(\pi) = \left(\bar{A}(\pi) - A(\pi) \right)+E_M(\pi)+B_M(\pi).$$ We analyze each of the three terms respectively. For the first term, under the assumption that $\mu(x,w)$ is Lipschitz on a bounded set, then $\sup_{x,w} \mu(x,w)$ is bounded. By Weak Law of Large Numbers \citep[Proposition 2.16]{vdV1998asymptotic}, $\bar{A}(\pi) - A(\pi) \stackrel{p}{\longrightarrow} 0$.

Let $\bar{\sigma}^2=\sup_{x,w} \sigma^2(x,w)$. It follows from  the Lipschitz assumption  that $\bar{\sigma}^2$ is bounded. Note that
\begin{equation*}
\begin{aligned}
 E\left[( \sqrt{n} E_M(\pi))^2\right] =& \frac{1}{n} E\left[\sum^n_{i=1} (2W_i-1)\left[(2\pi(X_i)-1)+\frac{K_M(\pi,i)}{M} \right]\epsilon_i\right]^2 \\ =& \nmeani E\left[\left(2\pi(X_i)-1 + \frac{\KMpi}{M}\right)^2 \epsilon_i^2 \right] \\ =& E\left[\left(2\pi(X_i)-1 + \frac{\KMpi}{M}\right)^2 \sigma^2(X_i,W_i) \right] \\ \leq& E\left[\left(1 + \frac{\KMi}{M}\right)^2 \sigma^2(X_i,W_i) \right] \\
 =& O(1),
\end{aligned}
\end{equation*}
where the second equality follows from that for $i \neq j$, $E\left[\epsilon_i \epsilon_j |\{X_i,W_i\}_{i=1}^n\right]=0$; the fourth inequality follows from triangle inequality and $|\KMpi| \leq \KMi$; and the last equality follows from the boundedness of $\bar{\sigma}^2$ and Lemma 3 in \citet{abadie2006large}. Hence, we obtain $ E_M(\pi) = O_p\left(n^{-1/2}\right)=o_p(1)$ by Markov's inequality. For the third term $B_M(\pi)$, it follows from Lemma \ref{order_conditional_bias} that $B_M(\pi)=O_p\left(n^{-1/p}\right)=o_p(1)$.

\end{proof}

\subsubsection{Proof of Theorem 1 (ii)}
\begin{proof}
Write $$\sqrt{n}\left(\hat{A}_{match}(\pi)-B_M(\pi)-A(\pi)\right)=\sqrt{n}\left(\bar{A}(\pi)-A(\pi)\right) + \sqrt{n} E_M(\pi).$$
We consider the asymptotic distributions for each term separately. For the first term, note that $$\left(2\pi(X)-1\right)\left(\mu(X,1)-\mu(X,0)\right)-A(\pi)$$ has mean zero and finite variance by the assumptions. Therefore, it follows from Central Limit Theorem \citep[Proposition 2.17]{vdV1998asymptotic} that
\begin{equation*}
\sqrt{n}\left(\bar{A}(\pi)-A(\pi)\right) \stackrel{d}{\rightarrow} N\left(0,V^{\tau(X)}(\pi)\right).
\end{equation*}

For the second term, we consider the similar arguments in the proof of \citet[Theorem 4]{abadie2006large}. Let $E_{M,i}(\pi) = (2W_i-1)\left[(2\pi(X_i)-1)+\frac{K_M(\pi,i)}{M} \right]\epsilon_i$. Conditional on $\{X_i,W_i\}_{i=1}^n$, $\KMpi$ is deterministic. Hence, the unit-level terms $\{E_{M,i}(\pi)\}^{n}_{i=1}$ are independently and not identically distributed with mean zero. We proceed to prove the asymptotic normality of $\sqrt{n} E_M(\pi) / \sqrt{V^E(\pi)}= \sum_{i=1}^n E_{M,i}(\pi)\sqrt{nV^E(\pi)}$ using Lindeberg-Feller Central Limit Theorem \citep[Theorem 27.2]{billingsley1995probability}. Specifically, given $\{X_i,W_i\}_{i=1}^n$, we need to verify the Lindeberg condition:
\begin{equation*}
	\frac{1}{nV^E(\pi)} \sum_{i=1}^n E\left[\left(E_{M,i}(\pi)\right)^2 I\left\{\left|E_{M,i}(\pi)\right|>\eta \sqrt{n V^E(\pi)}\right\} \mid \{X_i,W_i\}_{i=1}^n\right] \rightarrow 0, \text{ for all } \eta>0.
\end{equation*}
Note that
\begin{equation*}
\begin{aligned}
E& \left[\left(E_{M,i}(\pi)\right)^2 I\left\{\left|E_{M,i}(\pi)\right|>\eta \sqrt{n V^E(\pi)}\right\}|\{X_i,W_i\}_{i=1}^n\right] \\
& \leq E\left[\left(E_{M,i}(\pi)\right)^4\mid \{X_i,W_i\}_{i=1}^n\right]^{\frac{1}{2}} E\left[I\left\{\left|E_{M,i}(\pi)\right|>\eta \sqrt{n V^E(\pi)}\right\}\mid\{X_i,W_i\}_{i=1}^n\right]^{\frac{1}{2}} \\ & \leq E\left[\left(E_{M,i}(\pi)\right)^4\mid\{X_i,W_i\}_{i=1}^n\right]^{\frac{1}{2}} \Pr\left(\left|E_{M,i}(\pi)\right|>\eta \sqrt{n V^E(\pi)}\mid \{X_i,W_i\}_{i=1}^n\right) \\
& \leq  E\left[\left(E_{M,i}(\pi)\right)^4\mid \{X_i,W_i\}_{i=1}^n\right]^{\frac{1}{2}} \frac{E\left[\left(E_{M,i}(\pi)\right)^2\mid \{X_i,W_i\}_{i=1}^n\right]}{n\eta^2 V^E(\pi)},
\end{aligned}
\end{equation*}
where the first inequality follows from H\"older's inequality and the third inequality follows from Markov's inequality.

Let $\bar{c}=\sup_{x,w} E\left[\epsilon^4 | X=x, W=w \right] < \infty$ by assumption. Since $|\KMpi| \leq |\KMi|$, it follows that $|E_{M,i}(\pi)| \leq \left(1+\frac{\KMi}{M}\right) |\epsilon_i|$. In addition, we have $\underline{c}=\inf_{\textbf{X},\textbf{W}} V^E(\pi)>0$ by assumption. Therefore, we have
\begin{equation*}
\begin{aligned}
\frac{1}{nV^E(\pi)} & \sum_{i=1}^n E\left[\left(E_{M,i}(\pi)\right)^2 I\left\{\left|E_{M,i}(\pi)\right|>\eta \sqrt{n V^E(\pi)}\right\} \mid \{X_i,W_i\}_{i=1}^n\right] \\
& \leq \frac{1}{nV^E(\pi)} \sum_{i=1}^n \left[\left(1+\frac{\KMi}{M}\right)^4 E\left[\epsilon_i^4 \mid \{X_i,W_i\}_{i=1}^n \right]\right]^{\frac{1}{2}} \frac{\left(1+\frac{\KMi}{M}\right)^2 \sigma^2(X_i,W_i)}{n\eta^2 V^E(\pi)} \\
& \leq \frac{\bar{c}^{1/2} \bar{\sigma}^2}{\eta^2 {\underline{c}}^2}\frac{1}{n} \left[\nmeani \left(1+\frac{\KMi}{M}\right)^4\right].
\end{aligned}
\end{equation*}
It follows from the Lemma 3 of \citet{abadie2006large} that $E\left[\left(1+\frac{\KMi}{M}\right)^4\right]$ is bounded.  Therefore, we have
\begin{equation*}
\frac{1}{nV^E(\pi)} \sum_{i=1}^n E\left[\left(E_{M,i}(\pi)\right)^2 I\left\{\left|E_{M,i}(\pi)\right|>\eta \sqrt{n V^E(\pi)}\right\} \mid \{X_i,W_i\}_{i=1}^n\right] \stackrel{p}{\rightarrow} 0,
\end{equation*}
and thus $\sqrt{n} E_M(\pi) / \sqrt{V^E(\pi)} \stackrel{d}{\rightarrow} N(0,1)$.

Note that $\sqrt{n}\left(\bar{A}(\pi)-A(\pi)\right)$ and $\sqrt{n} E_M(\pi) / \sqrt{V^E(\pi)}$ are asymptotically independent. Using the same argument in \citet[Proof of Theorem 4, Page 266]{abadie2006large}, we complete the proof.
\end{proof}

\subsection{Proof of Theorem  2}
\begin{proof}
We only need to verify that $\sqrt{n}\left(B_M(\pi)-\hat{B}_M(\pi)\right) = o_p(1)$. To complete the proof, we need to use the result of unit-level bias correction \citep[Lemma A.2]{abadie2011bias}, which is given as follows:

\textit{Lemma A.2. \citet{abadie2011bias}}. Assume the conditions in Theorem 2  of \citet{abadie2011bias} hold, then
\begin{equation*}
\underset{i\leq n}{\max}  \left|\mu(X_i,w)-\mu(X_{j_m(i)},w)-\left(\hat{\mu}(X_i,w)-\hat{\mu}(X_{j_m(i)},w)\right)\right| = o_p\left(n^{-1/2}\right)
\end{equation*}
for $w=0,1$.

It follows from \eqref{BMpi_rewrite} that
\begin{equation*}
\begin{aligned}
	|B_M&(\pi)-\hat{B}_M(\pi)| \\
	=& \frac{1}{nM} \bigg| \sum_{i=1}^n \sum_{m=1}^M \left(2\pi(X_i)-1\right)\bigg\{W_i\bigg[\mu(X_i,0)-\mu(X_{j_m(i)},0)-\left(\hat{\mu}(X_i,0)-\hat{\mu}(X_{j_m(i)},0)\right)\bigg] \\ & - (1-W_i)\bigg[\mu(X_i,1)-\mu(X_{j_m(i)},1)-\left(\hat{\mu}(X_i,1)-\hat{\mu}(X_{j_m(i)},1)\right) \bigg]\bigg\} \bigg| \\
	\leq& \frac{1}{nM} \sum_{i=1}^n \sum_{m=1}^M \bigg[ \left|\mu(X_i,0)-\mu(X_{j_m(i)},0)-\left(\hat{\mu}(X_i,0)-\hat{\mu}(X_{j_m(i)},0)\right)\right| \\ & + \left|\mu(X_i,1)-\mu(X_{j_m(i)},1)-\left(\hat{\mu}(X_i,1)-\hat{\mu}(X_{j_m(i)},1)\right) \right| \bigg] \\
	\leq& \frac{1}{M} \sum_{m=1}^M \underset{i\leq n}{\max} \underset{w=0,1}{\sum} \left|\mu(X_i,w)-\mu(X_{j_m(i)},w)-\left(\hat{\mu}(X_i,w)-\hat{\mu}(X_{j_m(i)},w)\right)\right|.
\end{aligned}
\end{equation*}
Using Lemma A.2 in \citet{abadie2011bias}, we have $B_M(\pi)-\hat{B}_M(\pi)=o_p\left(n^{-1/2}\right)$.

\end{proof}

\subsection{Proof of Theorem 3}
\begin{proof}
	Since $R(\hat{\pi})=\left(A(\pi^{*})-A(\hat{\pi})\right)/2$ by definition, it suffices to analyze $A(\pi^{*})-A(\hat{\pi})$ throughout the proof.
	Let  $\Delta(\pi_1,\pi_2)=A(\pi_1)-A(\pi_2)$, $\bar{\Delta}(\pi_1,\pi_2)=\bar{A}(\pi_1)-\bar{A}(\pi_2)$, and $\hat{\Delta}(\pi_1,\pi_2)=\hat{A}^{bc}_{match}(\pi_1)-\hat{A}^{bc}_{match}(\pi_2)$. It follows from the triangle inequality that
	\begin{align}\label{three term}
		A(\pi^*)-A(\hat{\pi})=& \hat{A}^{bc}_{match}(\pi^*)-\hat{A}^{bc}_{match}(\hat{\pi})+A(\pi^*)-A(\hat{\pi})-\left(\hat{A}^{bc}_{match}(\pi^*)-\hat{A}^{bc}_{match}(\hat{\pi})\right)  \nonumber \\
		\leq& A(\pi^*)-A(\hat{\pi})-\left( \hat{A}^{bc}_{match}(\pi^*)-\hat{A}^{bc}_{match}(\hat{\pi})\right)  \nonumber \\
		=& \Delta(\pi^*,\hat{\pi}) - \hat{\Delta}(\pi^*,\hat{\pi}) \nonumber\\
		\leq& \underset{\pi_1,\pi_2 \in \Pi}{ \sup}|\Delta(\pi_1,\pi_2)-\hat{\Delta}(\pi_1,\pi_2)| \nonumber\\
		\leq& \underset{\pi_1,\pi_2 \in \Pi}{ \sup}|\bar{\Delta}(\pi_1,\pi_2)-\Delta(\pi_1,\pi_2)|+\underset{\pi_1,\pi_2 \in \Pi}{ \sup}|\bar{\Delta}(\pi_1,\pi_2)-\hat{\Delta}(\pi_1,\pi_2)| \nonumber\\
		\leq& \underset{\pi_1,\pi_2 \in \Pi}{ \sup}|\bar{\Delta}(\pi_1,\pi_2)-\Delta(\pi_1,\pi_2)|+ 2 \underset{\pi \in \Pi}{\sup}|B_M(\pi)-\hat{B}_M(\pi)| +2 \underset{\pi \in \Pi}{\sup}|E_M(\pi)|,
	\end{align}
	where the last inequality follows from
	$$\hat{\Delta}(\pi_1,\pi_2)-\bar{\Delta}(\pi_1,\pi_2)=B_M(\pi_1)-\hat{B}_M(\pi_1)+E_M(\pi_1)-[B_M(\pi_2)-\hat{B}_M(\pi_2)+E_M(\pi_2)].$$
	We will derive the non-asymptotic high probability bounds for the three terms in \eqref{three term} in Step 1-3 separately, and then combine these bounds to complete the proof in Step 4.
	
	\subsubsection{Step 1: High probability bound for the first term in \eqref{three term}}
	
	To derive this bound, we utilize Lemma 2 in \citet{zhouzhengyuan2023}. First, we rewrite $A(\pi)$ and $\bar{A}(\pi)$ in a vector-based form in order to apply their Lemma 2. Let $$
	\Gamma^{vec}_i=\left[
	\begin{aligned}
		\begin{array}
			{c}
			\mu(X_i,0) - \mu(X_i,1)\\
			\mu(X_i,1) - \mu(X_i,0)
		\end{array}
	\end{aligned}
	\right] \text{ and }
	\pi^{vec}(X_i)=\left[
	\begin{aligned}
		\begin{array}
			{c}
			1-\pi(X_i)\\
			\pi(X_i)
		\end{array}
	\end{aligned}
	\right], \text{ for } i=1,\ldots,n.
	$$
	Hence, we have $\bar{A}(\pi)=\nmeani \left<\Gamma^{vec}_i, \pi^{vec}(X_i)\right>$ and $A(\pi)=E\left[\left<\Gamma^{vec}_i, \pi^{vec}(X_i)\right>\right]$, where $\left<\cdot,\cdot\right>$ is the inner product between two vectors.
	Second, since $|Y|\leq L$ almost surely, the random variable $\mu(X,1)-\mu(X,0)$ is bounded almost surely, and so does $\Gamma_i^{vec}$. Third, according to the statement in \citet[Proof of Corollary 1]{zhouzhengyuan2023}, the entropy integral $\kappa(\Pi) \leq 2.5 \sqrt{\text{VC}(\Pi)}$.
	Finally, it follows directly from Lemma 2 in \citet{zhouzhengyuan2023} that, for any $\delta>0$, with probability at least $1-2\delta$, there exists constants $c_1, c_2>0$,
	\begin{equation}\label{zhou2023_bound}
	\underset{\pi_1,\pi_2 \in \Pi}{ \sup}|\bar{\Delta}(\pi_1,\pi_2)-\Delta(\pi_1,\pi_2)|\leq \left(c_1\sqrt{\text{VC}(\Pi)}+c_2+\sqrt{2\log \frac{1}{\delta}}\right)\sqrt{\frac{V^*}{n}}+o(\frac{1}{\sqrt{n}}),
	\end{equation} where
	$V^* =  4\sup_{\pi_1,\pi_2 \in \Pi} E\left[\left(\mu(X,1)-\mu(X,0)\right)^2 \left(\pi_1(X)-\pi_2(X)\right)^2 \right].$
	
	For the last term in \eqref{zhou2023_bound}, let $c_3>0$, then there exists $N(c_3)$ such that for all $n>N(c_3)$, the last term is upper bounded by $\frac{c_3}{\sqrt{n}}$. Therefore, for all $n>N(c_3)$, with probability at least $1-2\delta$, we have
	\begin{equation}\label{Step 1 bound}
		\underset{\pi_1,\pi_2 \in \Pi}{ \sup}|\bar{\Delta}(\pi_1,\pi_2)-\Delta(\pi_1,\pi_2)|\leq \left(c_1\sqrt{\text{VC}(\Pi)}+c_2+\sqrt{2\log \frac{1}{\delta}}\right)\sqrt{\frac{V^*}{n}}+\frac{c_3}{\sqrt{n}}.
	\end{equation}
	
	\subsubsection{Step 2: High probability bound for the second term in \eqref{three term}}
	
	It follows from the proof of Corollary \ref{bias-corrected estimator normality} that
	\begin{equation*}
	\sup_{\pi \in \Pi} |B_M(\pi)-\hat{B}_M(\pi)| \leq \frac{1}{M} \sum_{m=1}^M \underset{i\leq n}{\max} \underset{w=0,1}{\sum} \left|\mu(X_i,w)-\mu(X_{j_m(i)},w)-\left(\hat{\mu}(X_i,w)-\hat{\mu}(X_{j_m(i)},w)\right)\right|.
	\end{equation*}
    Using Lemma A.2 in \citet{abadie2011bias}, we have $\sup_{\pi \in \Pi}|B_M(\pi)-\hat{B}_M(\pi)|=o_p\left(\frac{1}{\sqrt{n}}\right)$. Let $c_4>0$, for any $\delta>0$, then there exists $N(c_4,\delta)$ such that for all $n>N(c_4,\delta)$, with probability at least $1-\delta$,
	\begin{equation}\label{Step 2 bound}
		2 \underset{\pi \in \Pi}{\sup}|B_M(\pi)-\hat{B}_M(\pi)| \leq \frac{c_4}{\sqrt{n}}.
	\end{equation}
	
	\subsubsection{Step 3: High probability bound for the third term in \eqref{three term}}
	
	By the definition of $E_M(\pi)$ and triangle inequality, we have
	\begin{align}\label{decompose to I-IV term}
		\underset{\pi \in \Pi}{\sup}|E_M(\pi)| \leq & \underset{\pi \in \Pi}{\sup}\left|\frac{1}{n}\sum^n_{i=1} (2W_i-1)(2\pi(X_i)-1)\epsilon_i\right| + \underset{\pi \in \Pi}{\sup}\left|\frac{1}{n}\sum^n_{i=1} (2W_i-1)\frac{K_M(\pi,i)}{M}\epsilon_i\right| \nonumber \\
		\leq & \left|\frac{1}{n}\sum^n_{i=1} (2W_i-1)\epsilon_i \right| +\left|\frac{1}{n}\sum^n_{i=1} (2W_i-1) \frac{K_M(i)}{M} \epsilon_i \right|+ 2 \underset{\pi \in \Pi}{\sup}\left|\frac{1}{n}\sum^n_{i=1} (2W_i-1)\pi(X_i)\epsilon_i\right| \nonumber\\ &+ \frac{2}{M}\underset{\pi \in \Pi}{\sup}\left|\frac{1}{n}\sum^n_{i=1} (2W_i-1)\left[\sum_{j=1}^n \pi(X_j)M_{ji}\right]\epsilon_i\right| \nonumber \\
		\equiv& \text{(I)+(II)+(III)+(IV)},
	\end{align}
	where the second inequality follows from the fact that
	$$K_M(\pi,i)=\sum^n_{j=1}\left(2\pi(X_j)-1\right)M_{ji}=
	2\sum_{j=1}^n\pi(X_j)M_{ji} - K_M(i).$$ We will derive high probability upper bounds for these four terms separately.
	
	\noindent \textbf{Step 3.1: Bound the (I) term in \eqref{decompose to I-IV term}}
	
	\noindent We bound (I) using Hoeffding's inequality \citep[Theorem 2.2.6]{RomanVershynin2018HDP}. Under the assumption that $|Y_i|\leq L$, we have $|\epsilon_i|=|Y_i - E[Y_i | X_i,W_i]| \leq 2L$ by the triangle inequality. In addition, note that $E[(2W_i-1)\epsilon_i]=E\left[(2W_i-1)E\left[\epsilon_i \mid X_i,W_i\right]\right]=0$. By Hoeffding's inequality, for any $t>0$,
	$$\Pr\left(\left|\frac{1}{n} \sum_{i=1}^n (2W_i-1)\epsilon_i\right| > \frac{t}{\sqrt{n}}\right) \leq  2\exp\left(\frac{-t^2}{8L^2}\right).$$ Let $\delta=\exp\left(\frac{-t^2}{8L^2}\right)$. We obtain
	\begin{equation}\label{I term bound}
		\Pr\left(\left|\frac{1}{n} \sum_{i=1}^n (2W_i-1)\epsilon_i\right| \leq 2\sqrt{2}L \sqrt{\frac{\log \frac{1}{\delta}}{n} }\right) \geq 1-2\delta.
	\end{equation}
	
	\noindent \textbf{Step 3.2: Bound the (II) term in \eqref{decompose to I-IV term}}
	
	\noindent For (II),  we also apply Hoeffding's inequality to derive the high probability bound. Given $\{X_i,W_i\}^n_{i=1}$, $K_M(i)$ is constant and $(2W_i-1)K_M(i)\epsilon_i$ is independent with mean zero. For any $t>0$, we have
	$$\Pr\left(\left|\frac{1}{n} \sum_{i=1}^n (2W_i-1)\frac{K_M(i)}{M}\epsilon_i \right|> \frac{t}{n} \sqrt{\sum_{i=1}^n K^2_M(i)} \mid \{X_i,W_i\}^n_{i=1}\right) \leq 2\exp\left(-\frac{M^2 t^2}{8L^2}\right).$$
	Taking expectation with respect to (w.r.t.) $\{X_i,W_i\}^n_{i=1}$, we  have $$ \Pr\left(\left|\frac{1}{n} \sum_{i=1}^n (2W_i-1)\frac{K_M(i)}{M}\epsilon_i\right|> \frac{t}{n} \sqrt{\sum_{i=1}^n K^2_M(i)}\right) \leq 2\exp\left(-\frac{M^2 t^2}{8L^2}\right).$$ Let $\delta=\exp\left(\frac{-M^2 t^2}{8L^2}\right)$. It follows that
	\begin{equation}\label{II term bound}
		\Pr\left(\left|\frac{1}{n} \sum_{i=1}^n (2W_i-1)\frac{K_M(i)}{M}\epsilon_i\right| \leq \frac{2\sqrt{2}L}{Mn} \sqrt{\log \frac{1}{\delta}\sum_{i=1}^n K^2_M(i)} \right) \geq 1-2\delta.
	\end{equation}
	
	\noindent \textbf{Step 3.3: Bound the (III) term in \eqref{decompose to I-IV term}}
	
	\noindent To derive this bound, we first analyze the order of $E \left[\sup_{\pi \in \Pi}\left|\frac{1}{n}\sum^n_{i=1} (2W_i-1)\pi(X_i)\epsilon_i\right|\right]$ using empirical process theory. We proceed to consider the empirical process indexed by the functional class $\Pi$. To simplify our analysis, we assume the cardinality of the class $\Pi$ is finite hereafter. Our results can be extended to a countable $\Pi$ by Monotone Convergence Theorem \citep{Massart2007concentration}.
	
	Denote the functional class $$\mathcal{H}=\left\{f_{\pi}(x,w,\epsilon)=(2w-1)\epsilon\pi(x), \pi \in \Pi\right\},$$ and the envelope function of $\mathcal{H}$ is $H(x,w,\epsilon)=|\epsilon|$. Therefore, we have $L_2(\mathbb{P})$ norm $\Vert H \Vert_{P,2}=\sqrt{E\epsilon^2} \leq 2L$. Recall that $E[(2W_i-1)\pi(X_i)\epsilon_i]=0$, we have
	\begin{equation*}
		\underset{\pi \in \Pi}{\sup}\left|\frac{1}{n}\sum^n_{i=1} (2W_i-1)\pi(X_i)\epsilon_i\right| = \underset{f \in \mathcal{H}}{\sup} |\mathbb{P}_nf|=\underset{f \in \mathcal{H}}{\sup} |(\mathbb{P}_n-\mathbb{P})f|,
	\end{equation*}
	where $\mathbb{P}_n$ and $\mathbb{P}$ denote the empirical measure and the expectation respectively.
	Note that the functional class $\mathcal{H}$ can be rewritten as the pointwise product between the class $\Pi$ and the function $g(w,\epsilon)=(2w-1)\epsilon$, i.e., $\mathcal{H}=\Pi \cdot g = \{\pi g, \pi \in \Pi\}$. Hence, by the results of VC-class \citep[Lemma 9.9]{kosorok2008introductionEP}, $\mathcal{H}$ is VC-class and $\text{VC}(\mathcal{H})\leq 2\text{VC}(\Pi)+1$ (note that Lemma 9.9 in \citet{kosorok2008introductionEP} is stated in terms of VC-index, which is equal to the corresponding VC dimension plus one).
	
	Since the norm of the empirical process indexed by a VC-class can be bounded by its VC dimension and the $L_2(\mathbb{P})$-seminorm of its envelope function \citep[Example 3.5.13]{vdV23}, then there exists a universal constant $c_5>0$ such that
	\begin{equation}\label{III term part 1}
		E\left[\underset{f \in \mathcal{H}}{\sup} |(\mathbb{P}_n-\mathbb{P})f|\right] \leq c_5 L \sqrt{\frac{\text{VC}(\Pi)}{n}}.
	\end{equation}
	Second, by Functional Hoeffding Theorem \citep[Theorem 3.26]{Wainwright2019HDS}, for any $t>0$ we have
	\begin{equation}\label{III term part 2}
		\begin{aligned}
		&\Pr\left(\underset{\pi \in \Pi}{\sup}\left|\frac{1}{n}\sum^n_{i=1} (2W_i-1)\pi(X_i)\epsilon_i\right| > E\left[\underset{\pi \in \Pi}{\sup}\left|\frac{1}{n}\sum^n_{i=1} (2W_i-1)\pi(X_i)\epsilon_i\right|\right]  + \frac{t}{\sqrt{n}} \right) \\
		&\leq \exp\left(-\frac{t^2}{64L^2}\right).
		\end{aligned}
	\end{equation}
	Let $\delta=\exp\left(-\frac{t^2}{64L^2}\right)$. Combining \eqref{III term part 1} and \eqref{III term part 2} we have
	\begin{equation}\label{III term bound}
		\Pr\left(\underset{\pi \in \Pi}{\sup}\left|\frac{1}{n}\sum^n_{i=1} (2W_i-1)\pi(X_i)\epsilon_i\right| \leq c_5 L \sqrt{\frac{\text{VC}(\Pi)}{n}}   + 8L \sqrt{\frac{\log \frac{1}{\delta}}{n}} \right) \geq 1-\delta.
	\end{equation}
	
	\noindent \textbf{Step 3.4: Bound the (IV) term in \eqref{decompose to I-IV term}}

	\noindent To do this, we first derive the order of
	\begin{equation}\label{IV term expectation only}
	E\left[	\underset{\pi \in \Pi}{\sup}\left|\frac{1}{n}\sum^n_{i=1} (2W_i-1)\left[\sum_{j=1}^n \pi(X_j)M_{ji}\right]\epsilon_i\right| \right]
	\end{equation}
    using Dudley’s entropy integral bound \citep[Theorem 8.1.3]{RomanVershynin2018HDP}. Then we consider the concentration inequality of $$\underset{\pi \in \Pi}{\sup}\left|\frac{1}{n}\sum^n_{i=1} (2W_i-1)\left[\sum_{j=1}^n \pi(X_j)M_{ji}\right]\epsilon_i\right|$$  around the expectation \eqref{IV term expectation only} using Talagrand's inequality \citep{Massart2007concentration}. In our proof, we need to use the maximal inequality of the absolute values of sub-Gaussian random variables, which is given as follows:

    \textit{ Corollary 7.2. \citet{chensongxiConcentrationReview}.} Let $\{X_i\}_{i=1}^n$ are sub-Gaussian random variables such that for $t \geq 0$, $\Pr(|X_i| > t) \leq 2\exp\left(-t^2/(2\sigma^2)\right)$,
    then
    \begin{equation}\label{maximal inequality subGaussian}
    E\left[\underset{ i \leq n}{\max} |X_i|\right] \leq \sigma \sqrt{2 \log(2n)}.
    \end{equation}
	
	\noindent \textbf{Step 3.4.1: The order of expectation}
	
	\noindent We define the unit-level term $a_i(\pi)$ in (IV) as $$ a_i(\pi)=(2W_i-1)\left[\sum_{j=1}^n \pi(X_j)M_{ji}\right]\epsilon_i, \text{ for } i=1,\ldots,n.$$ Given $\{X_i,W_i\}^n_{i=1}$, $a_i(\pi)$ is independently and not identically distributed with $E\left[a_i(\pi)\right]=0$. We define a Rademacher process $G^0_{n}(\pi)=\frac{1}{\sqrt{n}}\sum^n_{i=1} a_i(\pi) Z_i$, where $\{Z_i\}^n_{i=1}$ is a sequence of independent Rademacher variables  independent of $\{X_i,W_i,\epsilon_i\}^n_{i=1}$. Therefore, conditional on $\{X_i,W_i\}^n_{i=1}$, it follows from the symmetrization inequality \citep[Lemma 11.4]{Boucheron2013ConcentrationInequalities} that
	\begin{equation}\label{symmetrization}
		\begin{aligned}
			E\left[\underset{\pi \in \Pi}{\sup}\left|\frac{1}{n}\sum^n_{i=1} (2W_i-1)\left[\sum_{j=1}^n \pi(X_j)M_{ji}\right]\epsilon_i\right| \right] = & \frac{1}{\sqrt{n}} E\left[\underset{\pi \in \Pi}{\sup}\left|\frac{1}{\sqrt{n}}\sum^n_{i=1} a_i(\pi)\right| \right] \\
			\leq & \frac{2}{\sqrt{n}} E\left[\underset{\pi \in \Pi}{\sup}\left|\frac{1}{\sqrt{n}}\sum^n_{i=1} a_i(\pi) Z_i\right| \right] \\
			= &  \frac{2}{\sqrt{n}} E \left[\underset{\pi \in \Pi}{\sup} \left|G^0_n(\pi)\right|\right].
		\end{aligned}
	\end{equation}
	We proceed to derive an upper bound for the conditional expectation $E \left[\sup_{\pi \in \Pi} \left|G^0_n(\pi)\right|\right]$ with fixed $\{X_i,W_i\}^n_{i=1}$.
	
	For any function $f : \mathcal{X} \rightarrow \mathbb{R}$, the $L_2(\mathbb{P}_n)$-seminorm is defined as $\Vert f \Vert_n = \sqrt{\frac{1}{n} \sum_{i=1}^n f^2(X_i)}$ \citep{vdV23}.  Given $\{X_i,W_i,\epsilon_i\}^n_{i=1}$, the sequence $\{a_i(\pi)\}_{i=1}^n$ is deterministic. Therefore, by Hoeffding's inequality,  we have
	\begin{equation*}
		\Pr\left(\left| G^0_{n}(\pi_1) - G^0_{n}(\pi_2)\right| \geq t \right) \leq 2 \exp \left( -t^2 \bigg/ \left\{ \frac{2}{n} \sum_{i=1}^n \epsilon^2_i \left[
		\sum_{j=1}^n \left(\pi_1(X_j)-\pi_2(X_j)\right)M_{ji} \right]^2 \right\} \right),
	\end{equation*}
	for any $\pi_1, \pi_2 \in \Pi$ and any $t \geq 0$.
	Note that
	\begin{align*}
		\frac{1}{n} & \sum_{i=1}^n \epsilon^2_i \left[
		\sum_{j=1}^n \left(\pi_1(X_j)-\pi_2(X_j)\right)M_{ji} \right]^2
		= \frac{1}{n} \sum_{i=1}^n \epsilon^2_i \left[
		\underset{j: M_{ji}=1}{\sum}\left(\pi_1(X_j)-\pi_2(X_j)\right)\right]^2 \\
		\leq &
		\frac{1}{n} \sum_{i=1}^n \epsilon^2_i K_M(i) \left[
		\underset{j: M_{ji}=1}{\sum}\left(\pi_1(X_j)-\pi_2(X_j)\right)^2 \right] \\
		\leq & \frac{4L^2  \underset{i \leq n}{\max} K_M(i)}{n} \sum_{i=1}^n \left[
		\underset{j: M_{ji}=1}{\sum}\left(\pi_1(X_j)-\pi_2(X_j)\right)^2 \right] \\
		= & \frac{4L^2 \underset{i \leq n}{\max} K_M(i)}{n} \sum_{i=1}^n \left[
		\sum_{j=1}^n \left(\pi_1(X_j)-\pi_2(X_j)\right)^2 M_{ji} \right] \\
		= & \frac{4L^2 \underset{i \leq n}{\max} K_M(i)}{n} \sum_{j=1}^n \left[
		\left(\pi_1(X_j)-\pi_2(X_j)\right)^2  \sum_{i=1}^n M_{ji} \right] \\
		= & \frac{4ML^2  \underset{i \leq n}{\max} K_M(i)}{n} \sum_{j=1}^n \left[
		\left(\pi_1(X_j)-\pi_2(X_j)\right)^2  \right]
		=  4ML^2 \underset{i \leq n}{\max} K_M(i) \Vert \pi_1 - \pi_2 \Vert^2_n,
	\end{align*}
	where the first inequality follows from Cauchy-Swarchz inequality and the fact that $K_M(i)=\sum_{j=1}^n M_{ji}=\sum_{j: M_{ji}=1} 1$; the fourth equality follows from $\sum_{i=1}^n M_{ji}=M$, i.e., the cardinality of the set $\mathcal{J}_M(j)$ is $M$.
	Therefore, conditional on $\{X_i,W_i,\epsilon_i\}^n_{i=1}$,  $ G^0_{n}(\pi)$ is a  sub-Gaussian process with mean zero relative to $L_2(\mathbb{P}_n)$-seminorm $\Vert \cdot \Vert_n$ in the sense that
	\begin{equation}\label{sub-Gaussian G0 process}
		\Pr\left(\left| G^0_{n}(\pi_1) - G^0_{n}(\pi_2)\right| \geq t\right) \leq 2 \exp \left(\frac{-t^2}{8ML^2 \underset{i \leq n}{\max} K_M(i) \Vert \pi_1 - \pi_2 \Vert^2_n} \right),
	\end{equation}
    for any $\pi_1, \pi_2 \in \Pi$ and any $t \geq 0$.

	We will derive an upper bound for $E_Z \left[\sup_{\pi \in \Pi} \left|G^0_n(\pi)\right|\right]$ by slightly modifying the chaining technique in the proof of Dudley's entropy integral bound, where $E_Z(\cdot)$ denotes the expectation w.r.t. $Z$. We follow the similar argument in \citet{BSen2021empiricalprocess}, though the proof can be found elsewhere \citep[e.g.,][]{RomanVershynin2018HDP,Wainwright2019HDS,vdV23}.
	
	Let $\tilde{\pi}_0(\cdot) \equiv 0$, we consider the class $\tilde{\Pi}=\Pi \cup \{\tilde{\pi}_0\}$. The diameter of $\tilde{\Pi}$ is $\tilde{D}=\sup_{\pi^{\prime},\pi^{\prime \prime} \in \tilde{\Pi}} \Vert \pi^{\prime} -\pi^{\prime \prime} \Vert_n$, and note that $\tilde{D} \leq  1$. For $k \geq 1$, let $\tilde{\Pi}_k$ denote the maximal $\tilde{D}2^{-k}$-separated subset of $\tilde{\Pi}$ and hence the cardinality of $\tilde{\Pi}_k$ is $D(2^{-k},\tilde{\Pi},L_2(\mathbb{P}_n))$. By the maximality, it follows that
	\begin{equation}\label{maximality of subset}
	\sup_{\pi^{\prime} \in \tilde{\Pi}} \inf_{\pi^{\prime \prime} \in \tilde{\Pi}_k} \Vert \pi^{\prime} - \pi^{\prime \prime}\Vert_n\leq \tilde{D}2^{-k}.
	\end{equation}
    In addition, since $\tilde{\Pi}$ is finite under our assumption, there exists a sufficiently large $S$ such that $\tilde{\Pi}_S=\tilde{\Pi}$.  For any $\pi \in \tilde{\Pi}$, we can choose $\tilde{\pi}_k $ such that $\tilde{\pi}_k = \arg\inf_{\psi \in \tilde{\Pi}_k} \Vert \pi - \psi\Vert_n$, and therefore $\Vert \pi - \tilde{\pi}_k \Vert_n \leq \tilde{D}2^{-k}$ by  \eqref{maximality of subset}. Let  $\tilde{\Pi}_0=\{\tilde{\pi}_0\}$, we also have $\Vert \pi - \tilde{\pi}_0 \Vert_n \leq \tilde{D}$. Furthermore, we note that $\tilde{\pi}_S=\pi$ since  $\tilde{\Pi}_S=\tilde{\Pi}$. Finally, for any $\pi \in \tilde{\Pi}$, we can find a sequence of policies  $\tilde{\pi}_0,\tilde{\pi}_1,\ldots,\tilde{\pi}_S$ that approximates $\pi$.
	
	Note that $G^0_n(\tilde{\pi}_0)=0$ and that $G^0_n(\tilde{\pi}_S)=G^0_n(\pi)$, we rewrite the process $G^0_n(\pi)$ as $$G^0_n(\pi)=\sum_{k=1}^S \left(G^0_n(\tilde{\pi}_k)-G^0_n(\tilde{\pi}_{k-1})\right).$$
	It follows from the triangle inequality that
	\begin{equation}\label{chaining-expectation inequality}
		E_Z\left[\underset{\pi \in \Pi}{\sup} \left|G^0_n(\pi)\right|\right] = E_Z\left[\underset{\pi \in \tilde{\Pi}}{\sup} \left|G^0_n(\pi)\right|\right] \leq \sum_{k=1}^S E_Z\left[ \underset{\pi \in \tilde{\Pi}}{\sup} \left| G^0_n(\tilde{\pi}_k) - G^0_n(\tilde{\pi}_{k-1}) \right|\right].
	\end{equation}
	Recall that in \eqref{sub-Gaussian G0 process}, we have showed that $G^0_n(\pi)$ is a sub-Gaussian process given $\{X_i,W_i,\epsilon_i\}^n_{i=1}$. Noting that $\Vert \tilde{\pi}_k - \tilde{\pi}_{k-1} \Vert_n \leq \Vert \tilde{\pi}_k - \pi \Vert_n + \Vert \tilde{\pi}_{k-1} - \pi \Vert_n \leq 3\tilde{D}2^{-k}$, we have
	\begin{equation*}
	\Pr\left(\left| G^0_{n}(\tilde{\pi}_k) - G^0_{n}(\tilde{\pi}_{k-1})\right| \geq t\right)\leq 2 \exp \left(\frac{-t^2}{72 \cdot  2^{-2k} M\tilde{D}^2L^2 \underset{i \leq n}{\max} K_M(i)} \right).
	\end{equation*}
     Since the set of random variables $\left\{G^0_{n}(\pi^{\prime}) - G^0_{n}(\pi^{\prime \prime})\right\}_{\pi^{\prime} \in \tilde{\Pi}_k,\pi^{\prime \prime} \in \tilde{\Pi}_{k-1}}$ contains at most $$D(2^{-(k-1)},\tilde{\Pi},L_2(\mathbb{P}_n)) \cdot D(2^{-k},\tilde{\Pi},L_2(\mathbb{P}_n)) \leq  D^2(2^{-k},\tilde{\Pi},L_2(\mathbb{P}_n))$$ elements. By the maximal inequality of sub-Gaussian random variables \eqref{maximal inequality subGaussian}, we have
     \begin{equation*}
     E_Z\left[ \underset{\pi \in \tilde{\Pi}}{\sup} \left| G^0_n(\tilde{\pi}_k) - G^0_n(\tilde{\pi}_{k-1}) \right|\right] \leq 2^{-k} C \tilde{D} L \sqrt{M\underset{i \leq n}{\max} K_M(i)} \sqrt{\log\left(2D(2^{-k},\tilde{\Pi},L_2(\mathbb{P}_n))\right)},
     \end{equation*}
	 and therefore by \eqref{chaining-expectation inequality} we have (recall that the generic constant $C$ may varies from place to place)
	 \begin{equation*}
	 \begin{aligned}
	 E_Z\left[\underset{\pi \in \tilde{\Pi}}{\sup} \left|G^0_n(\pi)\right|\right] \leq &
	 \sum_{k=1}^S 2^{-k} C\tilde{D} L \sqrt{M\underset{i \leq n}{\max} K_M(i)} \sqrt{\log\left(2D(2^{-k},\tilde{\Pi},L_2(\mathbb{P}_n))\right)} \\
	 = &2 C  L \sqrt{M\underset{i \leq n}{\max} K_M(i)} \sum_{k=1}^S  \int_{\tilde{D}2^{-k-1}}^{\tilde{D}2^{-k}}  \sqrt{\log\left(2D(2^{-k},\tilde{\Pi},L_2(\mathbb{P}_n))\right)} d \eta \\
	 \leq & C \tilde{D} L \sqrt{M\underset{i \leq n}{\max} K_M(i)} \sum_{k=1}^S  \int_{\tilde{D}2^{-k-1}}^{\tilde{D}2^{-k}}  \sqrt{\log\left(2D(\eta,\tilde{\Pi},L_2(\mathbb{P}_n))\right)} d \eta \\
	 \leq & C L \sqrt{M\underset{i \leq n}{\max} K_M(i)}  \int_{0}^{\tilde{D}/2}  \sqrt{\log\left(2D(\eta,\tilde{\Pi},L_2(\mathbb{P}_n))\right)} d \eta \\
	 \leq & C L \sqrt{M\underset{i \leq n}{\max} K_M(i)}  \int_{0}^{\tilde{D}/2}  \sqrt{\log D(\eta,\tilde{\Pi},L_2(\mathbb{P}_n))} d \eta \\
	 \leq & C L \sqrt{M\underset{i \leq n}{\max} K_M(i)}  \int_{0}^{\tilde{D}/2}  \sqrt{\log N(\eta/2,\tilde{\Pi},L_2(\mathbb{P}_n)) } d \eta \\
	 \leq & C L \sqrt{M\underset{i \leq n}{\max} K_M(i)}  \int_{0}^{1/4}  \sqrt{\log N(\eta,\tilde{\Pi},L_2(\mathbb{P}_n)) } d \eta,
	 \end{aligned}
	 \end{equation*}
	 where the fourth inequality follows from the fact that, for $\eta \leq \tilde{D}/2$, the packing number $D(\eta,\tilde{\Pi},L_2(\mathbb{P}_n)) \geq 2$; the fifth inequality follows from the relationship between covering number and packing number \citep[Page 147]{vdV23}; and the last inequality follows from the change of variable and the fact that $\tilde{D} \leq 1$.
	
	Let  $Q$ denote all the finitely discrete probability measures $Q$ on $\mathcal{X}$, note that $$ N(\eta,\tilde{\Pi},L_2(\mathbb{P}_n)) \leq 1 + N(\eta,\Pi,L_2(\mathbb{P}_n)),$$ and by the result of covering number of VC-class functions \citep[Theorem 2.6.7]{vdV23}, we have
	\begin{equation}\label{entropy-intergal}
	\begin{aligned}
	E_Z \left[\underset{\pi \in \tilde{\Pi}}{\sup} \left|G^0_n(\pi)\right|\right] \leq &  CL \sqrt{M\underset{i \leq n}{\max} K_M(i)}  \int_{0}^{1/4}  \sqrt{1+\log N(\eta,\Pi,L_2(\mathbb{P}_n)) } d \eta \\
	\leq & CL\sqrt{M \underset{i \leq n}{\max} K_M(i)} \int_{0}^{1/4} \sqrt{1+\underset{Q}{\sup}\log N(\eta,\Pi,L_2(Q))} d \eta \\
	\leq & CL\sqrt{M \underset{i \leq n}{\max} K_M(i)  \text{VC}(\Pi)}.
	\end{aligned}
	\end{equation}
	Taking expectation w.r.t. $\{X_i,W_i,\epsilon_i\}^n_{i=1}$ in \eqref{entropy-intergal}, by  Jensen's inequality we have
	\begin{equation}\label{entropy-intergal-jensen}
		E \left[\underset{\pi \in \Pi}{\sup} \left|G^0_n(\pi)\right| \right] \leq  E \left[CL\sqrt{M  \underset{i \leq n}{\max} K_M(i)  \text{VC}(\Pi)}\right]
		\leq CL \sqrt{M \cdot \text{VC}(\Pi) E\left[\underset{i \leq n}{\max} K_M(i)\right]}.
	\end{equation}
	
	Since $E[K^q_M(i)]$ is bounded uniformly in $n$ for any $q>0$ \citep[Lemma 3]{abadie2006large}, and $\{K_M(i)\}^n_{i=1}$ is identically distributed, then by the sharper maximal inequality \citep[Corollary 7.1]{chensongxiConcentrationReview}, for any $q > 1$ we have
	\begin{equation}\label{maximal-KMi}
		E\left[\underset{i \leq n}{\max} K_M(i)\right] = o(n^{\frac{1}{q}}).
	\end{equation}
	Therefore, for some $c_6>0$, there exists $N(c_6,q)$ such that for all $n>N(c_6,q)$, we have $E\left[\max_{i \leq n} K_M(i)\right] \leq c_6 n^{\frac{1}{q}}$.
	Combining \eqref{symmetrization}, \eqref{entropy-intergal-jensen} and \eqref{maximal-KMi}, for any $q>1$ and $n>N(c_6,q)$, we have
	\begin{equation}\label{IV term expectation}
		\begin{aligned}
			E\left[\underset{\pi \in \Pi}{\sup}\left|\frac{1}{n}\sum^n_{i=1} (2W_i-1)\left[\sum_{j=1}^n \pi(X_j)M_{ji}\right]\epsilon_i\right|\right] \lesssim & L \sqrt{\frac{M \cdot \text{VC}(\Pi)}{n^{1-1/q}}},
		\end{aligned}
	\end{equation}
	where we use the notation $A \lesssim B$ to denote that $A$ is less than $B$ up to a universal constant.
	
	\noindent \textbf{Step 3.4.2: Concentration inequality}
	
	\noindent Let $b_i(\pi)=(2W_i-1)\left[\sum_{j=1}^n \pi(X_j)M_{ji}\right]$. Conditional on $\{X_i,W_i\}_{i=1}^n$, the term $$\sum \limits_{i=1}^n (2W_i-1)\left[\sum_{j=1}^n \pi(X_j)M_{ji}\right] \epsilon_i = \sum \limits_{i=1}^n b_i(\pi)\epsilon_i $$ is a sum of independent terms and  $E\left[b_i(\pi)\epsilon_i\right]=0$. Note that
	\begin{equation*}
		\left|\frac{b_i(\pi)\epsilon_i}{2L \cdot \underset{i \leq n}{\max} K_M(i)}\right| \leq  \frac{\sum_{j=1}^n \pi(X_j)M_{ji}}{ \underset{i \leq n}{\max} K_M(i)} \leq \frac{\sum_{j=1}^n M_{ji}}{\underset{i \leq n}{\max} K_M(i)}=\frac{K_M(i)}{\underset{i \leq n}{\max} K_M(i)} \leq 1.
	\end{equation*}
	Given $\{X_i,W_i\}_{i=1}^n$, the unit-level terms $\left\{\frac{b_i(\pi)\epsilon_i}{2L \cdot \max_{i \leq n} K_M(i)}\right\}_{i=1}^n$ satisfy the conditions of Talagrand's inequality \citep[Page 169-170]{Massart2007concentration}. Therefore, it follows that for any $ t >0$, with at least probability $1-\exp(-t)$, we have
	\begin{equation*}
		\begin{aligned}
			\underset{\pi \in \Pi}{\sup}  \left|\sum_{i=1}^n\frac{b_i(\pi)\epsilon_i}{2L \cdot \underset{i \leq n}{\max} K_M(i)}\right| \leq &  E_{\epsilon}\left[ \underset{\pi \in \Pi}{\sup}  \left|\sum_{i=1}^n\frac{b_i(\pi)\epsilon_i}{2L \cdot \underset{i \leq n}{\max} K_M(i)}\right| \right] + 2t \\
			& +2\sqrt{\left(2V+16E_{\epsilon}\left[ \underset{\pi \in \Pi}{\sup}  \left|\sum_{i=1}^n\frac{b_i(\pi)\epsilon_i}{2L \cdot \underset{i \leq n}{\max} K_M(i)}\right| \right]\right) t},
		\end{aligned}
	\end{equation*}
	where
	\begin{equation*}
		V=\underset{\pi \in \Pi}{\sup} \sum_{i=1}^n E_{\epsilon}\left[\frac{b^2_i(\pi) \epsilon^2_i}{4L^2 \cdot \left[ \underset{i \leq n}{\max} K_M(i)\right]^2 } \right] \leq \underset{\pi \in \Pi}{\sup} \frac{\sum_{i=1}^n b^2_i(\pi)}{\left[ \underset{i \leq n}{\max} K_M(i)\right]^2} \leq \frac{\sum_{i=1}^n K_M^2(i)}{\left[ \underset{i \leq n}{\max} K_M(i)\right]^2},
	\end{equation*}
	and $E_\epsilon(\cdot)$ denote the expectation of $\epsilon$. Note that for any $\alpha >0$,
	\begin{equation*}
		\begin{aligned}
		&	2\sqrt{\left(2V+16E_{\epsilon}\left[ \underset{\pi \in \Pi}{\sup}  \left|\sum_{i=1}^n\frac{b_i(\pi)\epsilon_i}{2L \cdot \underset{i \leq n}{\max} K_M(i)}\right| \right]\right) t}\\
        \leq & 2\sqrt{2V t} + 8 \sqrt{E_{\epsilon}\left[ \underset{\pi \in \Pi}{\sup}  \left|\sum_{i=1}^n\frac{b_i(\pi)\epsilon_i}{2L \cdot \underset{i \leq n}{\max} K_M(i)}\right| \right] t} \\
			\leq & 2\sqrt{2V t} + \alpha E_{\epsilon}\left[ \underset{\pi \in \Pi}{\sup}  \left|\sum_{i=1}^n\frac{b_i(\pi)\epsilon_i}{2L \cdot \underset{i \leq n}{\max} K_M(i)}\right| \right] + \frac{16t}{\alpha}.
		\end{aligned}
	\end{equation*}
	Without loss of generality, we set $\alpha=1$. It follows that, conditional on $\{X_i,W_i\}^n_{i=1}$, the following inequality holds with probability  at least $1-\exp(-t)$:
	\begin{equation*}
            \underset{\pi \in \Pi}{\sup}  \left|\frac{1}{n} \sum_{i=1}^n  b_i(\pi)\epsilon_i\right| \leq
			2E_{\epsilon}\left[\underset{\pi \in \Pi}{\sup}  \left|\frac{1}{n} \sum_{i=1}^n  b_i(\pi)\epsilon_i\right|\right] +  \frac{36Lt}{n}\underset{i \leq n}{\max} K_M(i) +  \frac{4L\sqrt{2t}}{n}\sqrt{\sum_{i=1}^n K^2_M(i)}.
	\end{equation*}

	Since $\{b_i(\pi)\}_{i=1}^n$ are constants given $\{X_i,W_i\}^n_{i=1}$, we replace the $E_{\epsilon}\left[\sup_{\pi \in \Pi}  \left|\frac{1}{n} \sum_{i=1}^n  b_i(\pi)\epsilon_i\right|\right]$ above by $E\left[\sup_{\pi \in \Pi}  \left|\frac{1}{n} \sum_{i=1}^n  b_i(\pi)\epsilon_i\right|\right]$. Taking expectation w.r.t. $\{X_i,W_i\}^n_{i=1}$, for any $t>0$, with  probability at least $1-\exp(-t)$, we have
	\begin{equation}\label{IV term Talagrand concentration}
		\underset{\pi \in \Pi}{\sup}  \left|\frac{1}{n} \sum_{i=1}^n  b_i(\pi)\epsilon_i\right| \leq
		2E\left[\underset{\pi \in \Pi}{\sup}  \left|\frac{1}{n} \sum_{i=1}^n  b_i(\pi)\epsilon_i\right|\right] + \frac{36Lt}{n}\underset{i \leq n}{\max} K_M(i) +  \frac{4L\sqrt{2t}}{n}\sqrt{\sum_{i=1}^n K^2_M(i)}.
	\end{equation}
	Let $\delta=\exp(-t)$. It follows from \eqref{IV term expectation} and \eqref{IV term Talagrand concentration} that for all $n > N(c_6,q)$, with  probability at least $1-\delta$, we have
	\begin{equation}\label{IV term final bound}
		\underset{\pi \in \Pi}{\sup}  \left|\frac{1}{n} \sum_{i=1}^n  b_i(\pi)\epsilon_i\right| \lesssim
		L \sqrt{\frac{M \cdot \text{VC}(\Pi)}{n^{1-1/q}}} + \frac{L \log\frac{1}{\delta}}{n}\underset{i \leq n}{\max} K_M(i) +  \frac{L\sqrt{\log\frac{1}{\delta}}}{n}\sqrt{\sum_{i=1}^n K^2_M(i)} .
	\end{equation}
	
	\subsubsection{Step 4: High probability bound for $A(\pi^*)-A(\hat{\pi})$}
	
	Combining \eqref{Step 1 bound}, \eqref{Step 2 bound}, \eqref{decompose to I-IV term}, \eqref{I term bound}, \eqref{II term bound}, \eqref{III term bound} and \eqref{IV term final bound}, it follows that for any $\delta>0$, with probability at least $1-9\delta$ and for all $n> \max\{N(c_3), N(c_4,\delta), N(c_6,q)\}$, we have
	\begin{equation}\label{Total bound 1}
		\begin{aligned}
			A(\pi^*)-A(\hat{\pi}) \lesssim & \frac{L}{\sqrt{M}}\sqrt{\frac{\text{VC}(\Pi)}{n^{1-1/q}}} + \frac{1}{\sqrt{n}}\left(\sqrt{\log\frac{1}{\delta}}+\sqrt{\text{VC}(\Pi)}\right) \left(\sqrt{V^*}  + L\right)  \\ &+ \frac{L\sqrt{\log\frac{1}{\delta}}}{Mn}\sqrt{\sum_{i=1}^n K^2_M(i)} + \frac{L\log\frac{1}{\delta}}{Mn}\underset{i \leq n}{\max} K_M(i).
		\end{aligned}
	\end{equation}
    Note that $\frac{1}{n}\sqrt{\sum_{i=1}^n K^2_M(i)}=\frac{1}{\sqrt{n}} \sqrt{\frac{1}{n}\sum_{i=1}^n K^2_M(i)} = O_p\left(\frac{1}{\sqrt{n}}\right)=o_p\left(\frac{1}{\sqrt{n^{1-1/q}}}\right)$, where the second equality follows from  Markov's inequality. For $c_7>0$, then for any $\delta>0$, there exists $N(c_7,\delta)$ such that for all $n>N(c_7,\delta)$, we have
	\begin{equation}\label{Total bound 2}
		\Pr\left(	\frac{1}{n}\sqrt{\sum_{i=1}^n K^2_M(i)} \leq \frac{c_7}{\sqrt{n^{1-1/q}}}\right) \geq 1-\delta.
	\end{equation}
	Similarly, since $\frac{1}{n}\max_{i \leq n} K_M(i)=o_p\left(\frac{1}{n^{1-1/q}}\right) $ by Markov's inequality. For  $c_8>0$, there exists $N(c_8,\delta)$ such that for all $n>N(c_8,\delta)$, we have
	\begin{equation}\label{Total bound 3}
		\Pr\left(	\frac{1}{n}\underset{i \leq n}{\max} K_M(i) \leq \frac{c_8}{n^{1-1/q}}\right) \geq 1-\delta.
	\end{equation}
	
	Recall that $A(\pi^{*})-A(\hat{\pi})=2R(\hat{\pi})$. Let $N_{\delta,q}= \max\{N(c_3), N(c_4,\delta), N(c_6,q),N(c_7,\delta),N(c_8,\delta)\}$. It follows from \eqref{Total bound 1}, \eqref{Total bound 2}, and \eqref{Total bound 3}  that, for any $\delta>0$,  $q>1$ and $n > N_{\delta,q}$, with  probability at least $1-11\delta$,  the following inequality holds:
	\begin{equation*}
	R(\hat{\pi}) \lesssim
	\frac{L}{\sqrt{M n^{1-1/q}}} \left(\sqrt{\text{VC}(\Pi)} + \sqrt{\frac{\log\frac{1}{\delta}}{M}}\right) +\frac{1}{\sqrt{n}}\left(\sqrt{\text{VC}(\Pi)}+\sqrt{\log\frac{1}{\delta}}\right)\left(\sqrt{V^*}  + L\right) + \frac{L \log\frac{1}{\delta}}{M n^{1-1/q}}.
\end{equation*}
\end{proof}

\newpage
\section{Additional Information for Simulation Studies}
\subsection{Data Generating Process}
Table \ref{data_generating_table} presents detailed information for the data generating processes in simulation studies. Additionally, the proportions of the controlled versus treated units are displayed. For the propensity score model,  Scenario 1 is the same as \citet{kangschafer2007}, where $l(\cdot)$ is a linear combination of predictors  and the proportion is balanced. In Scenario 2 we consider a nonlinear form of $l(\cdot)$, while the proportion is still balanced. Scenario 3 is modified from Scenario 1 to introduce imbalanced proportion. We consider a randomized experiment in Scenario 4, where the proportion is extremely imbalanced. Scenario 5 is modified from \citet{wu2020matched}, where a nonlinear propensity score model is in use and the  proportions in the two groups are  extremely imbalanced.
\begin{table}[H]
	\centering
	\caption{Data generating processes for different scenarios in simulation studies. The outcome model is $Y = m(X) + Wc(X) + e$, where the binary treatment variable $W$ is generated by $\Pr(W=1|X)=\exp\left(l(X)\right)/\left[1+\exp\left(l(X)\right)\right]$.}
	\label{data_generating_table}
	\resizebox{0.8\textwidth}{!}{\begin{tabular}{llc}
			\toprule
			Propensity score model & \makecell[c]{$l(X)$} & Controlled versus treated \\
			\midrule
			1. Linear, balanced & $-X_1+0.5X_2-0.25X_3-0.1X_4$ &  50:50 \\
			2. Nonlinear, balanced & $0.1X_1^3+0.2X_2^3+0.3X_3$  & 50:50   \\
			3. Linear, imbalanced & $2.1-X_1+2X_2-0.25X_3-0.1X_4$ & 23:77  \\
			4. Constant, imbalanced & $\log 9$ & 10:90 \\
			5. Nonlinear, imbalanced & $1+\exp(X_2)+\sin X_1 \cos X_3$ & 12:88 \\
			\midrule
			Main effect model & \makecell[c]{$m(X)$} &   \\
			\midrule
			1. Linear & $1 + 2 X_1 - X_2 + 0.5 X_3 - 1.5 X_4 $ & \\
			2. Nonlinear & $4\sin(X_1) + 2.5\cos(X_2) - X_3X_4$ & \\
			\midrule
			Contrast function & \makecell[c]{$c(X)$} &  \\
			\midrule
			1. Tree  & $2I\left\{X_1>0, X_2>0\right\}-1$ & \\
			2. Non-tree & $2I\left\{2X_2-\exp(1+X_1)+2 >0\right\}-1$ & \\
			\bottomrule
	\end{tabular}}
\end{table}

\newpage
\subsection{Simulation Results with Nonlinear Main Effect}
\begin{figure}[H]
	\centering
	\includegraphics[width=1\linewidth]{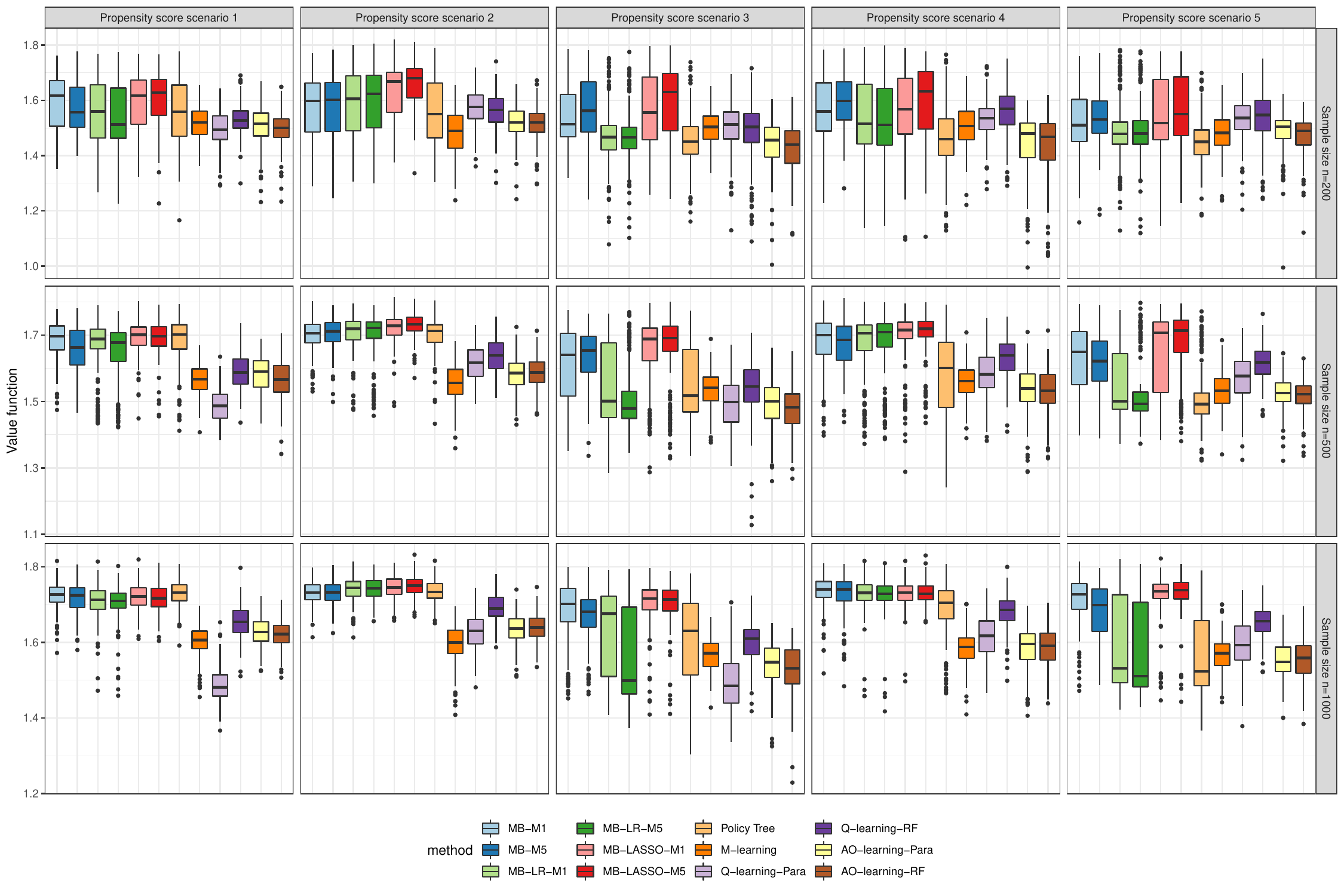}
	\caption{Boxplot of empirical value functions in different scenarios, where the main effect model is nonlinear and the contrast function is tree. The global optimal value function is 1.77.}
	\label{performance-nonlinearmain-tree}
\end{figure}

\begin{figure}[H]
	\centering
	\includegraphics[width=1\linewidth]{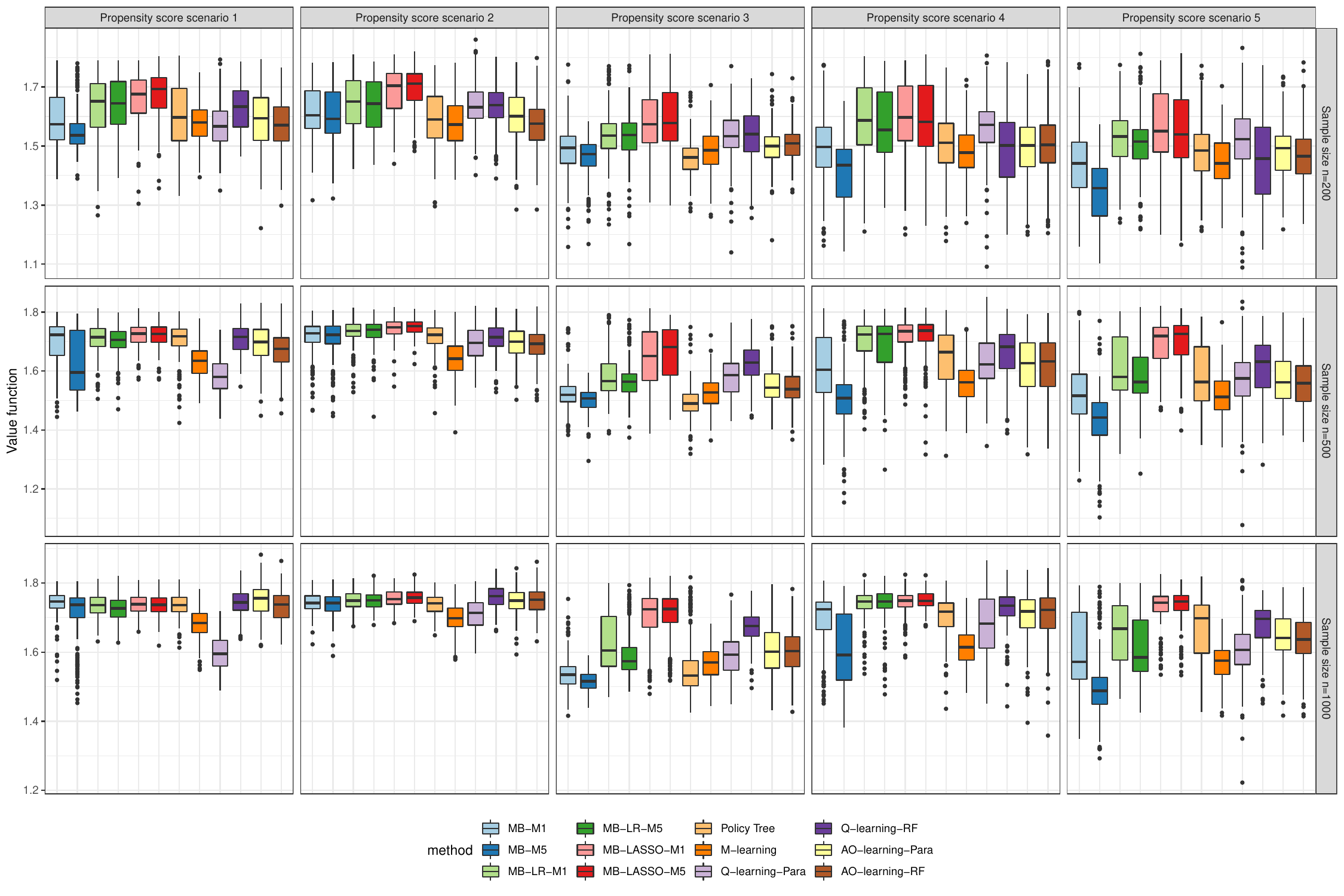}
	\caption{Boxplot of empirical value functions in different scenarios, where the main effect model is nonlinear and the contrast function is non-tree. The global optimal value function is 1.87.}
	\label{performance-nonlinearmain-nottree}
\end{figure}

\newpage
\subsection{Comparison of MB-Learning and Policy Tree with Sample Size Varying from 200 to 8000}
\begin{figure}[H]
	\centering
	\includegraphics[width=1\linewidth]{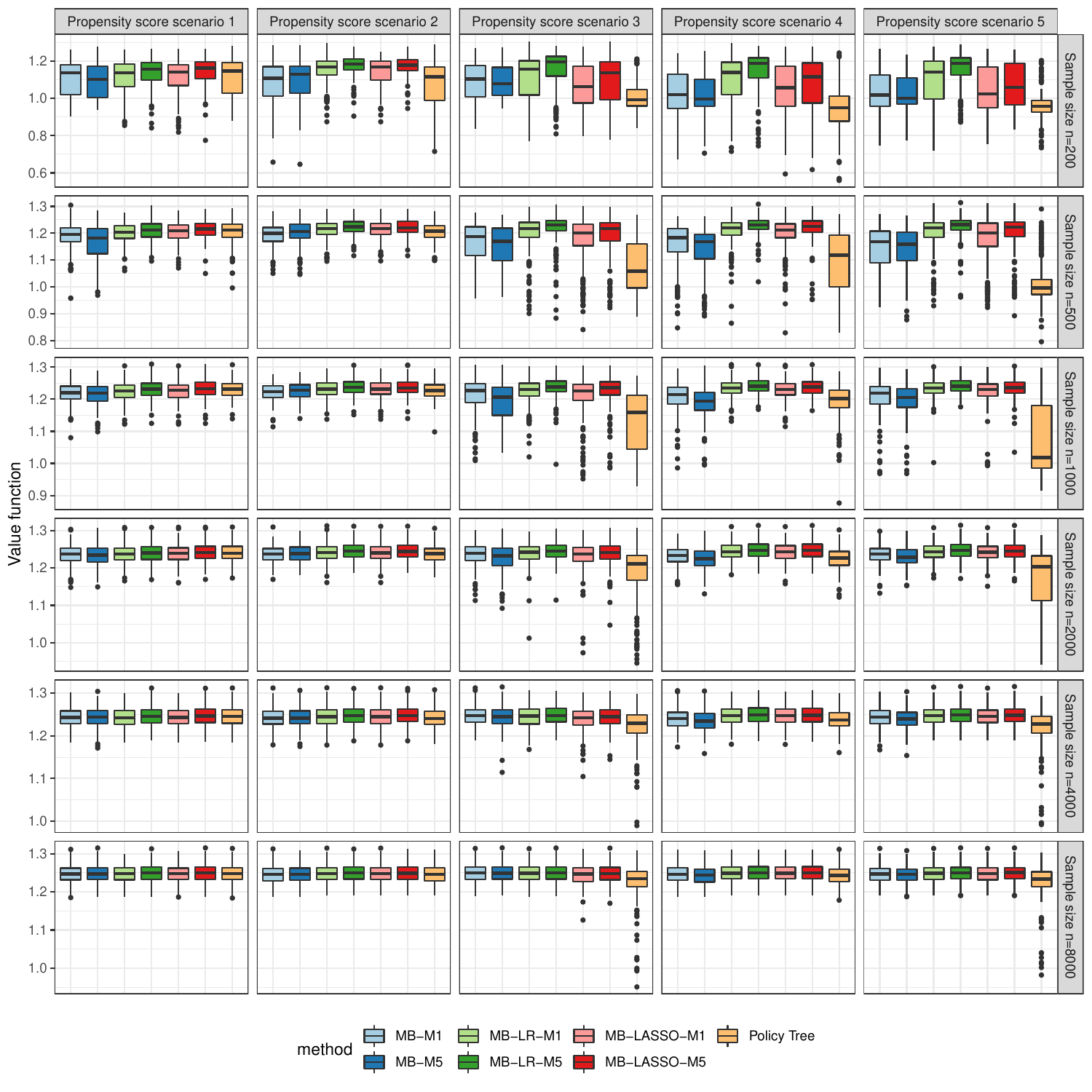}
	\caption{Boxplot of empirical value functions in different scenarios, where the main effect model is linear and the contrast function is tree.}
	\label{compare-with-policytree-trueboundary-tree}
\end{figure}

\begin{figure}[H]
	\centering
	\includegraphics[width=1\linewidth]{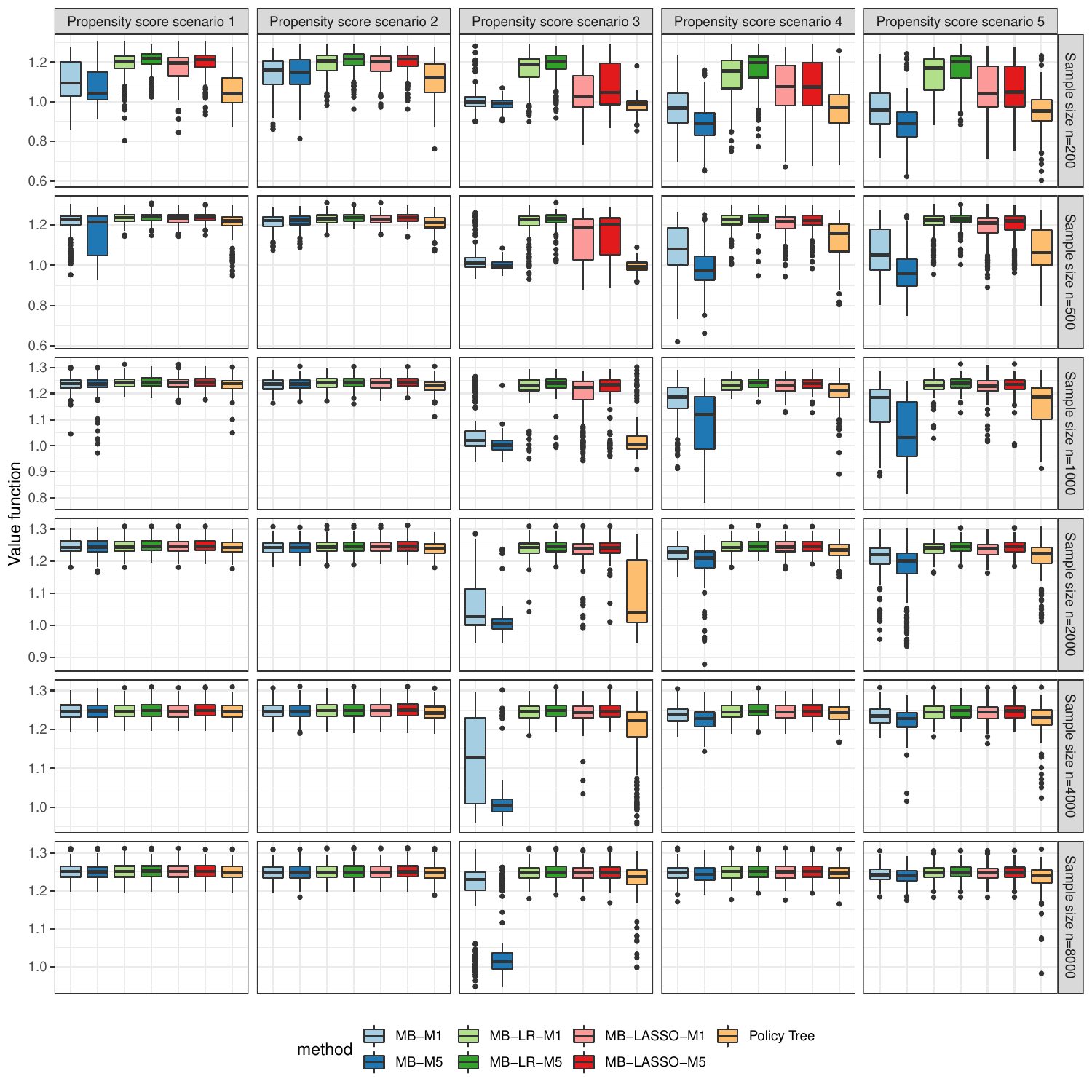}
	\caption{Boxplot of empirical value functions in different scenarios, where the main effect model is linear and the contrast function is non-tree.}
	\label{compare-with-policytree-trueboundary-nottree}
\end{figure}

\newpage
\section{Additional Results in Real Data Application}
\begin{figure}[H]
	\centering
	\begin{subfigure}{0.45\textwidth}
	\centering
	\includegraphics[width=1\linewidth]{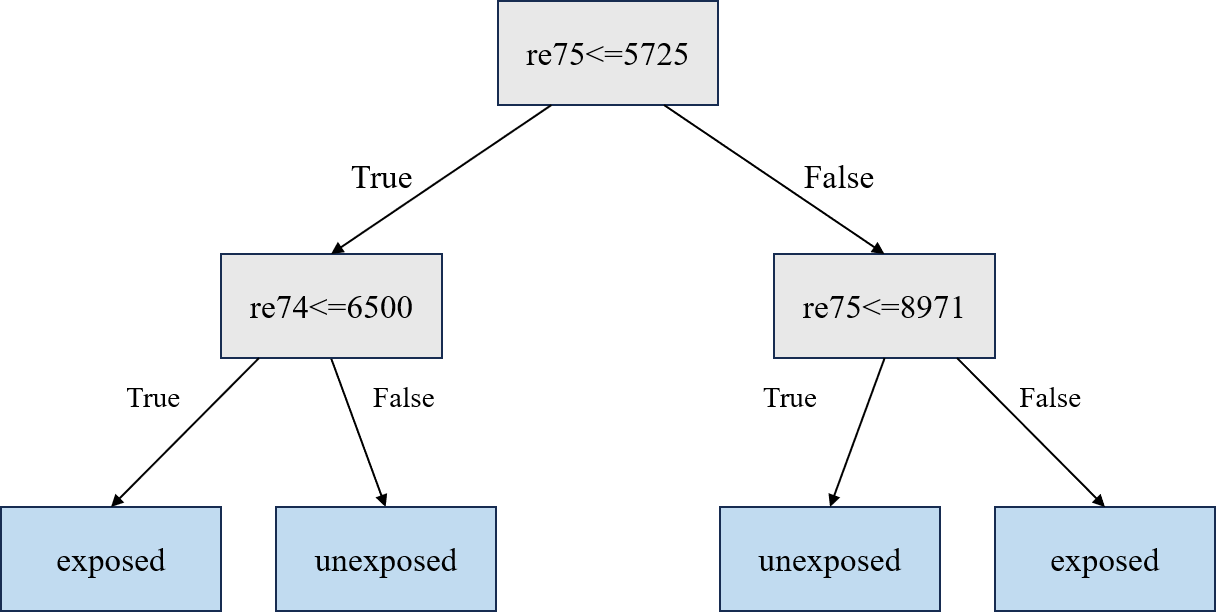}
	\caption{Depth-2 decision tree learned by MB-LASSO-M5}
	\end{subfigure}
    \hspace{5mm}
	\begin{subfigure}{0.45\textwidth}
	\centering
	\includegraphics[width=1\linewidth]{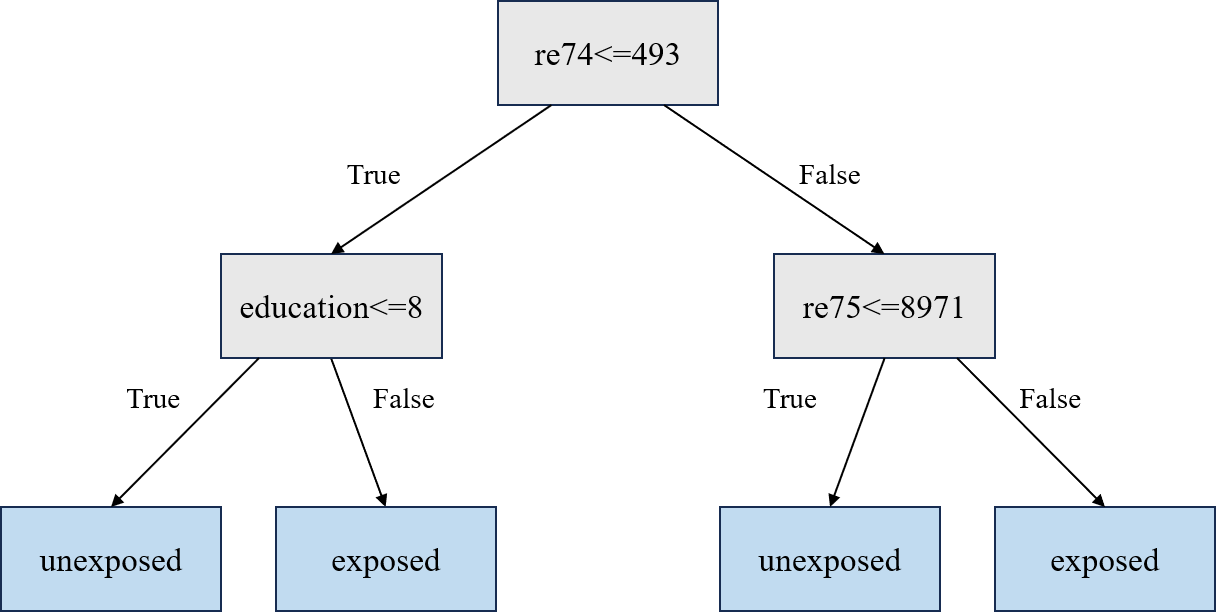}
	\caption{Depth-2 decision tree learned by Policy Tree}
    \end{subfigure}
    \caption{Optimal policies estimated by using all the data in DW dataset. (a) The optimal policy derived from MB-LASSO-M5 depends only on pre-treatment earnings, and 401 (90\%) of the job-seekers are recommended to be "exposed". (b) The optimal policy derived from Policy Tree depends on both pre-treatment earnings and education, and 290 (65\%) of the  job-seekers are recommended to be "exposed".}
    \label{learned-policy-DW}
\end{figure}

\vspace{1em}

\begin{figure}[H]
	\centering
	\begin{subfigure}{0.45\textwidth}
		\centering
		\includegraphics[width=1\linewidth]{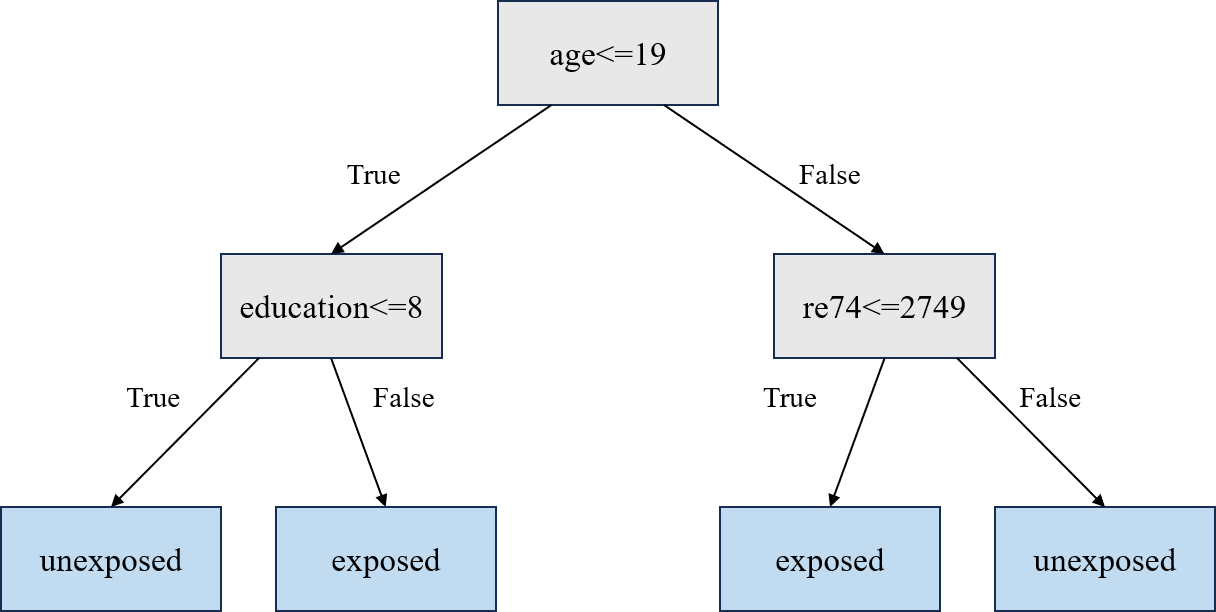}
		\caption{Depth-2 decision tree learned by MB-M5}
	\end{subfigure}
    \hspace{5mm}
	\begin{subfigure}{0.45\textwidth}
		\centering
		\includegraphics[width=1\linewidth]{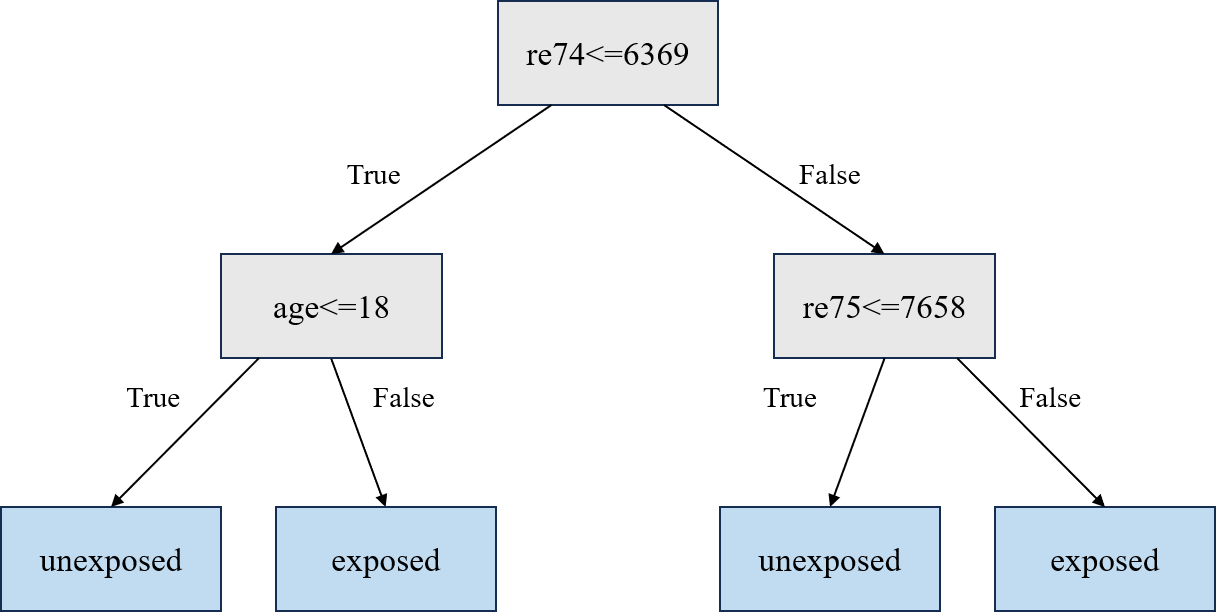}
		\caption{Depth-2 decision tree learned by Policy Tree}
	\end{subfigure}
	\caption{Optimal policies estimated by using all the data in  DW-CPS3 dataset. (a) The optimal policy derived from MB-M5 depends on  age, education and pre-treatment earnings, and 560 (64\%) of the samples are recommended to be "exposed". (b) The optimal policy derived from Policy Tree depends on age and pre-treatment earnings, and 557 (64\%) of the samples are recommended to be "exposed".}
	\label{learned-policy-DW-CPS3}
\end{figure}

\end{document}